\documentclass[11pt]{article}

\usepackage{latexsym}
\usepackage[all]{xy}
\usepackage{amsmath}
\usepackage{amssymb}
\usepackage{amsfonts}
\usepackage{accents}

\newcommand {\h} {\underline{Hom}}
\newcommand {\OO} {\mathcal{O}}
\newcommand {\T} {\mathcal{T}}
\newcommand {\D} {\mathcal{D}}
\newcommand {\Top} {\mathcal{S}}

\newcommand  {\dgcat}     {dg-\mathbf{Cat}}
\newcommand  {\mdgcat}     {dg-\mathbf{Cat}^{\otimes}}
\newcommand {\scat} {\infty-\mathbf{Cat}}
\newcommand  {\cdga}     {\mathbf{cdga}}
\newcommand  {\dga}     {\mathbf{dga}}
\newcommand  {\dgl}     {\mathbf{dgl}}

\newcommand  {\dAff}     {\mathbf{dAff}}
\newcommand  {\ncdga}     {\mathbf{cdga}^{\leq 0}}

\newcommand  {\dSt}   {\mathbf{dSt}}
\newcommand {\s}{\infty}
\newcommand {\Opcat}{\mathcal{O}p-\mathcal{C}at}

\newtheorem{thm}{Theorem}[section]
\newtheorem{prop}[thm]{Proposition}
\newtheorem{lem}[thm]{Lemma}
\newtheorem{slem}[thm]{Sub-lemma}
\newtheorem{df}[thm]{Definition}
\newtheorem{cor}[thm]{Corollary}

\begin{document}

\title{\textbf{Operations on derived moduli spaces of branes}}  
\author{\bigskip\\
Bertrand To\"en\thanks{This work is partially supported by the ANR grant ANR-09-BLAN-0151 (HODAG).}
 \\
\small{I3M, Universit\'e de Montpellier 2}\\
\small{Case Courrier 051}\\
\small{Place Eug\`ene Bataillon}\\
\small{34095 Montpellier Cedex, France}\\
\small{e-mail: btoen@math.univ-montp2.fr}}

\date{June 2013}

\maketitle

\begin{abstract}
The main theme of this work is the study of the operations that naturally exist 
on moduli spaces of maps $Map(S,X)$, also called the space of branes of $X$ with respect $S$.
These operations will be constructed as operations on the (quasi-coherent) derived
category $\D(Map(S,X))$, in the particular case where $S$ has some close relations
with an operad $\OO$. 
More precisely, for an $\s$-operad $\OO$ and an algebraic variety $X$ (or more generally
a derived algebraic stack), satisfying some natural conditions, we prove that 
$\OO$ acts on the object $\OO(2)$ by mean cospans. This universal action is used to prove that 
$\OO$ acts on the derived category of the space of maps $Map(\OO(2),X)$, which will call
the brane operations. We apply the existence of these operations, as well as
their naturality in $\OO$, in order to propose a sketch for a proof of the \emph{higher formality conjecture}, 
a far reaching extension of Konstevich's formality's theorem. By doing so we present a positive answer
to a conjecture of Kapustin (see \cite[p. 14]{kap}), relating polyvector fields on a 
variety $X$ and deformations of 
the mono\"\i dal derived category $\D(X)$.  
\end{abstract}

\tableofcontents

\section*{Introduction}

This work is an investigation of the operations that exist naturally on the space of \emph{branes} from
the point of view of derived categories of sheaves. The \emph{branes}
of the title are by definition morphisms $S \longrightarrow X$, where $S$ is a sphere (or some of its
generalisations) and $X$ is an algebraic variety (or some of its generalisations). 
The space of all branes, for fixed domain $S$, form a moduli space $Map(S,X)$, and the purpose of 
this work is to formalise the natural operations that can be defined on its
quasi-coherent derived category $\D(Map(S,X))$. For obvious
reasons automorphisms and more generally endomorphisms of $S$ do act on the 
space $Map(S,X)$ and thus on the category $\D(Map(S,X))$ by functoriality. We will be mainly interested
in operations on $\D(Map(S,X))$ which are not induced by endomorphisms of $S$ but rather come from 
\emph{cospans} on $S$. The cospans on $S$ are diagrams of the form 
$\xymatrix{S \ar[r] & S' & \ar[l] S}$ which
produce correspondences at the level of moduli spaces
$$\xymatrix{Map(S,X) &  Map(S',X) \ar[r]^-{p}  \ar[l]_-{f} & Map(S,X),}$$
and thus provide endofunctors at the level of derived categories
$$p_{*}f^{*} : \D(Map(S,X)) \longrightarrow \D(Map(S,X)).$$

The content of the present work is two folds and proposes 
solutions to the following two general problems.

\begin{enumerate}

\item Formalize the general fact that cospans on $S$ define natural operations on 
$\D(Map(S,X))$.

\item Provide a systematic way to construct interesting instances of 
operations on $\D(Map(S,X))$ out of simple input data. 

\end{enumerate}

The solution to the above two problems presented in this work is heavily based
on the one side on the notion of $\s$-operads of \cite{cm1,cm2,mw,ha} and on the other side
on techniques of $\s$-topos and of derived algebraic geometry in the sense of \cite{seat,hagI,hagII,ht}.
In a nutshell our main result states that for an algebraic variety $X$, or more generally
a derived Artin stack, and for an $\s$-operad $\OO$ satisfying some simple conditions, 
the derived category $\D(Map(\OO(2),X))$ carries a canonical action of $\OO$. Applied to 
the special cases where $\OO$ is the $E_k$-operad and its degenerate version the
$P_k$-operad, we deduce from this theorem a higher formality theorem, relating 
iterated Hochschild cohomology of $X$ to shifted polyvector fields on $X$, a result that has
been conjectured in the introduction of \cite{ptvv}. The special case of $k=2$, 
together with the relations between higher Hochschild cohomology and the deformation theory
of iterated mono\"\i dal categories (see \cite{fr}), provides
a positive answer to a conjecture of A. Kapustin appearing in his 2010 ICM lecture (see \cite[p. 14]{kap}). \\

\begin{center}\textbf{Cospans and $\s$-operads}\end{center}

The formalism of cospans on $S$ can be developed inside any general $\s$-topos $\T$. For
an object $S \in \T$, cospans are diagrams in $\T$ of the form
$\xymatrix{S \ar[r] & S' & \ar[l] S.}$
These diagrams are organised into several structures. First of all they form 
an $\s$-category, where morphisms are commutative diagrams in $\T$
$$\xymatrix{S \ar[r] \ar[d]_-{id} & S' \ar[d] & \ar[l] S \ar[d]^-{id} \\
S \ar[r] & S'' &  S. \ar[l]}$$
They can also be composed by push-outs over $S$, and thus form 
a mono\"\i dal $\s$-category. The cospans also possess $n$-ary 
versions by considering more general diagrams of the form
$$\xymatrix{\displaystyle{\coprod_{n}} S \ar[r] & S' & \ar[l] S,}$$
for which now the composition is operadic and defines a \emph{categorical $\s$-operad} (see definitions \ref{d1} and \ref{d1'}).
Finally, cospans can form family over a base object $x\in \T$, as diagrams
$\xymatrix{S\times x \ar[r] & S' & \ar[l] S \times x,}$
in the comma $\s$-topos $\T/x$. 

The structure that subsumes these operations all together
 is that of a \emph{categorical $\s$-operad in $\T$} (see definition \ref{d1'}), 
which technically speaking consists of an $\s$-functor
$$\Omega^{op} \times \T^{op} \longrightarrow \scat,$$
where $\Omega$ is a certain category of rooted trees, and satisfying some exactness conditions
in the first variable and a descent condition in the second variable.
An important first construction of this work is that of a cospans operad $\underline{End}^{\sqcup}(S)$, 
which is a categorical $\s$-operad in $\T$ encoding all the possible operations on cospans over a fixed
object $S\in \T$ (see definition \ref{d3}). 

Our first result is that many $\s$-operads $\OO$ in $\T$ act on the object $\OO(2)$ by means of 
cospans, a fact which is expressed by the following theorem. 

\begin{thm}{(See \ref{t1})}\label{ti1}
Let $\OO$ be an $\s$-operad in $\T$ satisfying the following two conditions.
\begin{enumerate}
\item $\OO(0)\simeq \OO(1)\simeq *$

\item For all pair of integers $n,m\geq 2$, the following square
$$\xymatrix{
\OO(n+1) \times \OO(m) \displaystyle{\coprod_{\OO(2) \times \OO(n) \times \OO(m)}} 
\OO(n) \times \OO(m+1) \ar[r] \ar[d] & \OO(n+m) \ar[d] \\
\OO(n) \times \OO(m) \ar[r] & \OO(n+m-1)}$$
is cartesian in the $\s$-topos $\T$.
\end{enumerate}
Then there exists a natural morphism of categorical $\s$-operads in $\T$
$$ad : \OO \longrightarrow \underline{End}^{\sqcup}(\OO(2)),$$
or in other words $\OO$ acts naturally on $\OO(2)$ by cospans. For an integer $n$, the morphism 
$\OO(n) \longrightarrow \underline{End}^{\sqcup}(\OO(2))(n)$ corresponds to the following family of $n$-ary 
cospan on $\OO(2)$ over $\OO(n)$
$$\xymatrix{
\displaystyle{\coprod_{n}}\OO(n) \times \OO(2) \ar[r] \ar[dr] & \OO(n+1) \ar[d] & \ar[l] \OO(2) \times \OO(n) 
\ar[ld] \\
 & \OO(n). & }$$
\end{thm}

The construction of the morphism $ad$ is by no means easy and occupies the first half of the present work. It should be noted
that the conditions of the theorem are satisfied by many interesting $\s$-operads, as for instance the $E_{n}$-operads
as well as their variations (see \S 6). Another source of examples is given by modular operads of pointed stable curves
of genus zero. 
For these two examples the morphism $ad$ of the theorem above is easy to identified: for the  $E_{n}$-operad $ad$ is
an incarnation of the fact that 
$E_{n}$ acts on the $(n-1)$-dimensional shere $S^{n-1}$ by means of cobordisms. The case of the modular operad
is even more clear: a pointed stable curve is a cospans on the point in the category of schemes. A striking consequence
of the theorem \ref{ti1} is that these two facts are special cases of a very general phenomenon, and in particular
are purely operadic data (e.g. the $E_{n}$-operad already contains the information about its actions on $S^{n-1}$).

\begin{center} \textbf{Brane operations} \end{center}

The theorem \ref{ti1} will  be used in order to construct operations on the derived category of branes. 
We still work in a general $\s$-topos $\T$, and we assume that we are given 
a notion of derived categories as an $\s$-functor
$$\D : \T^{op} \longrightarrow \mdgcat_{pr},$$
from $\T^{op}$ into the $\s$-category of presentable symmetric mono\"\i dal dg-categories (see \S 5). With these data, 
and for all nice enough dg-category $A$, we define a categorical $\s$-operad $\underline{End}^{\otimes}(A)$
in $\T$ whose $n$-ary operations are the dg-functors
$$A^{\otimes n} \longrightarrow A,$$
as well as their deformations over objects of $\T$ (see definition \ref{d8} for the precise definitions). The theorem \ref{ti1}, together with 
the fact that all cospans in dg-categories represents morphisms (because every morphism has an adjoint by definition), 
leads to the following corollary, which is our second main result.

\begin{cor}{(See \ref{cp4})}\label{ci1}
Let $\T$ and $\D$ be as above, and let $X \in \T$ and $\OO$ be a categorical $\s$-operad
satisfying the conditions of the theorem \ref{ti1}. Under some technical finiteness conditions, there exists a natural morphism 
of categorical $\s$-operad in $\T$
$$B : \OO \longrightarrow  \underline{End}^{\otimes}(\D(B_{\OO}(X))),$$
where $B_{\OO}(X)$ is the space of branes of $\OO$ with respect to $X$ defined as the hom object
$$B_{\OO}(X):=Map(\OO(2),X) \in \T.$$
\end{cor}

The morphism $B$ is called the \emph{brane operations morphism} and encodes the action of $\OO$ on 
$\D(B_{\OO}(X))$ induced from its actions on $\OO(2)$ by our theorem \ref{ti1}. In more concrete terms, 
the morphism $B$ provides for all $n$ a dg-functor
$$\psi^{n} : \D(B_{\OO}(X))^{\otimes n} \longrightarrow \D(B_{\OO}(X)) \otimes \D(\OO(n)),$$
which should be thought as a \emph{family of multilinear dg-endofunctors of $\D(B_{\OO}(X))$ parametrized
by $\OO(n)$}. The various dg-functors $\psi^{n}$ are related by higher coherences expressing compatibilities
between compositions in $\OO$ and compositions of dg-functors. Tracking back the construction of the
morphism $ad$ we can explicitly described the dg-functor $\psi^{n}$ as follows. For a fixed $n$, 
there is a family of cospans over $\OO(n)$ given by the diagram
$$\xymatrix{
\displaystyle{\coprod_{n}} \OO(n) \times \OO(2) \ar[r] \ar[rd] & \OO(n+1) \ar[d] & \ar[l] \ar[ld] \OO(2) \times \OO(n) \\
 & \OO(n), & }$$
defined by the various compositions in $\OO$. Applying the internal hom construction 
in $\T/\OO(n)$ provides a family of spans over $\OO(n)$
$$\xymatrix{
 B_{\OO}(X)^{n}\times \OO(n) \ar[rd] & \ar[l]_-{q} Map_{/\OO(n)}(\OO(n+1),X) \ar[d] \ar[r]^-{p} & \ar[ld] B_{\OO}(X) \times \OO(n) \\
 & \OO(n). & }$$
The dg-functor $\psi^{n}$ is obtained as push-forward pull-back 
$$\psi^{n}:=p_{*}q^{*} : \D(B_{\OO}(X))^{\otimes n} \longrightarrow \D(B_{\OO}(X)) \otimes \D(\OO(n)).$$
The technical conditions mentioned in the corollary insure that the dg-functor $p_{*}$ and $q^{*}$ are nice
enough in order to behave well (e.g. the base-change formula is required at some point). \\

\begin{center}\textbf{Some consequences} \end{center}

In the last part of this work we present one possible application of our results, theorem \ref{ti1} and
its corollary \ref{ci1}. 

For this we specify the situation to the case where $\T$ is the
$\s$-topos of derived stacks over the site of derived affine schemes (over a field
of caracteristic zero) with the \'etale topology (see \cite{seat,hagII}). 
The object $X$ will be taken to be a derived Artin $1$-stack of the form $[Y/G]$, where
$Y$ is a quasi-projective derived scheme and $G$ is an algebraic group. The $\s$-operad
$\OO$ will be taken either to be the $E_{k}$-operad or its degenerate
vesion the $\mathbb{P}_{k}$-operad classifying ($(1-k)$-shifted) Poisson dg-algebras. By formality, the
$\s$-operads $E_{k}$ and $\mathbb{P}_{k}$ are essentially equivalent when $k>1$, and thus 
so are their brane operations. In our last section we explain how this
equivalence at the level of brane operations can be used in order to deduce the following formality theorem
(see \S 6).

\begin{cor}{(See \ref{cform})}\label{ci2}
Let $X$ be a derived Artin stack as above and $k>1$. Then, a choice of a formality equivalence 
$E_{k}\simeq \mathbb{P}_{k}$, 
determines a quasi-isomorphism of dg-Lie algebras
$$HH^{E_{k}}(X)[k] \simeq \mathbb{R}\Gamma(X,Sym_{\OO_{X}}(\mathbb{T}_{X}[-k]))[k].$$
\end{cor}

In this corollary $HH^{E_{k}}(X)$ denotes the $E_{k}$-Hochschild cohomology of $X$, defined 
as the self-ext of the iterated diagonals of $X$ (see also \cite{fr}). By the so-called Deligne's conjecture, which 
is now a theorem due to \cite{ha,fr}, $HH^{E_{k}}(X)$ is an $E_{k+1}$-algebra and thus
$HH^{E_{k}}(X)[k]$ possesses a natural dg-Lie algebra structure. The right hand side is a derived version 
of (shifted) polyvector fields on $X$, and is endowed with a dg-Lie algebra structure 
induced from the bracket of vector fields. This corollary is an extension to higher Hochschild
cohomology of the famous Kontsevich's formality theorem, which would be the case $k=1$ and $X$ a smooth 
variety. Using the well known relations between $HH^{E_{k}}(X)[k]$ and the deformation theory 
of $E_{k}$-algebras (see for example \cite{fr}), this corollary offers a positive answer to a
conjecture of Kapusting in \cite[p. 14]{kap} (take $k=2$ and $X$ a smooth variety). In general this corollary is 
a solution to the higher formality conjecture mentioned in the introduction of \cite{ptvv}, at least for $n\neq 
0$. It has important consequences, as for example the existence of deformation quantization of shifted
Poisson structures, which according to the results of \cite{ptvv} will lead to 
many instances of deformations of $E_{k}$-mono\"\i dal dg-categories of sheaves on 
derived moduli spaces. \\

\begin{center} \textbf{Relations to other works} \end{center}

A major motivation for this work has been the notion of \emph{categorified Gromov-Witten theory}, 
which is a rather old idea consisting of lifting the Gromov-Witten invariants from cohomology 
to derived categories (this idea is for instance mentioned in \cite{seat}, and has
also been independantly advertised recently by Y. Manin \cite{man}). One of the main problem in trying to 
define categorified GW invariants comes from the existence of higher homotopies and
the precise way in which they should be organised. The language of the present work 
clarifies this point: categorified GW invariants for an algebraic variety $X$
essentially consists of a morphism of $\s$-operads
in derived stacks
$$\{\bar{\mathcal{M}}_{g,n+1}\}_{g,n} \longrightarrow \underline{End}^{\otimes}(\D(X)),$$
 where $\{\bar{\mathcal{M}}_{g,n+1}\}_{g,n}$ is the $\s$-operad
(in Deligne-Mumford stacks) of pointed stable curves. 
However, things are not quite as simple and several modifications of the notions and constructions of the
present work should be performed in order to formally construct a morphism of $\s$-operads as above. 
First of all the modular operad $\{\bar{\mathcal{M}}_{g,n+1}\}$ is certainly of configuration type in the sense 
close to our definition \ref{d6} but does not satisfy our crucial hypothesis $\OO(0)=\OO(1)=*$, except in 
genus zero. Moreover, some of the finiteness conditions of our corollary \ref{cp4}
do not hold, due to the non quasi-compactness of the stacks of all stable maps to a fixed target.
Finally, we did not inlcude any stability conditions in our construction, and therefore 
our results deal with the whole mapping spaces and do not 
take into accound the nicer subspaces of stable objects
as this is necessary for GW theory.
Nevertheless, we do beleive that some non trivial modifications of our  constructions and results
will provide a way to construct the categorified Gromov-Witten invariants. A work in progress by 
E. Mann and M. Robalo aims to extend the techniques of this paper from $\s$-operads
to the more general setting of cyclic $\s$-operads and modular $\s$-operads, 
as well as to include stability conditions in the constructions, and will probably lead to a full 
construction of the categorified Gromov-Witten invariants in  a near future \\

The notion of brane operations on $\D(Map(\OO(2),X))$, as well 
its related notion of $\OO$-cohomology $HH^{\OO}(X)$ described in our \S 6, seem very closely related
to the notions of \emph{centers} in the sense of \cite{bfn} and \cite{ha}. The precise relations 
for a general $\s$-operad has not been investigated, one of the reason being that 
these two notions does make sense in slightly different settings. However, when $\OO=E_{k}$, 
the $E_{k}$-mono\"\i dal structure on $\D(Map(S^{k-1},X))$ deduced from our 
results is most certainly the same as the one on the iterated center of
$\D(X)$ studied in \cite{bfn}, even thought we have not tried to write down a precise comparison. 

There are many situations in which the brane operations obtained from our corollary \ref{ci1}
can be identified with some previously known notions or constructions. An example
appears in the derived geometric Satake correspondence, 
for which the derived dg-category of certain constructible sheaves on 
the affine Grassmanian of a reductive group $G$ is identified with the dg-category $\D(Map(S^{2},B(^{L}G))$
where $^{L}G$ is the Langlands dual group.
The natural mono\"\i dal structure induced from convolution on the affine Grassmanian 
provides an $E_{3}$-mono\"\i dal structure on $\D(Map(S^{2},B(^{L}G))$, which is precisely the
one given by the corollary \ref{ci1} ($\OO=E_{3}$, $X=B(^{L}G)$) (this example has been 
explained to me by J. Lurie almost ten years ago).

Finally, the results of this work are very much related to the work in progress with T. Pantev, 
M. Vaqui\'e and G. Vezzosi, about shifted symplectic and Poisson structures as well 
as deformation quantization in the derived setting. The application of the techniques of the present 
work to formality and deformation quantization will be further investigated in great details
in a near future (see \cite{ptvv2}).\\

\begin{center} \textbf{Organisation of the present work} \end{center}

We have organised this work into five sections. The first section does not pretend to bring something new, 
and contains some reminders and
preliminaries concerning $\s$-operads as well as categorical operads. 
We remind the formalism of $\s$-operads in the sense 
of \cite{cm1,cm2,mw}, defined as certain functors on trees, and put this in the 
$\s$-categorical framework. In a section paragraph we define the notion of
a categorical operad, which is a slight generalisation of that of a 
(unital and associative) mono\"\i dal category (we could have called 
these \emph{reduced multi-2-categories}). We explain how categorical operads are
related to $\s$-operads by means of a rectification construction, analog to the
construction which produce a strict mono\"\i dal category out of a mono\"\i dal 
category. 

In the second section we present a detailed construction of the $\s$-operads of cospans
$\underline{End}^{\sqcup}(X)$, of spans $\underline{End}^{\times}(Z)$, as well 
as the duality morphism 
$$\underline{End}^{\sqcup}(X) \longrightarrow \underline{End}^{\times}(Z^{X}).$$
Here the objects $X$ and $Z$ belongs to a given $\s$-topos $\T$, and these $\s$-operads
are $\s$-operads inside the $\s$-category of stacks of $\s$-categories over $\T$ (i.e. they
are family of $\s$-operads in $\s$-categories parametrized by objects in $\T$). These objects
are foundamental for our main result \ref{ti1}, and surprisingly already their definitions
require some non-trivial work. Our constructions are heavily based for instance on the notions 
of fibered $\s$-categories (see e.g. \cite{tove2}).

The third section concernes the proof of our first theorem \ref{ti1}. We start by 
a lax version (see theorem \ref{t1}), which states that in absence of condition $(2)$ 
the morphism $ad$ exists as a lax morphisms of categorical $\s$-operads. The condition $(2)$ 
is then shown to be equivalent to the fact that the lax morphism $ad$ is an actual
morphism. The proof is presented in 
several intermediate steps and is again heavily based on the notions and techniques of
fibered $\s$-categories. The theorem \ref{t1} is also used to produce 
a dual version of the morphism $ad$. 

In the section 4 we explain how the morphism $ad$ of theorem \ref{ti1} produces 
operations on the derived category $\D(Map(\OO(2),X))$, and thus a proof of
corollary \ref{ci1}. For this we 
assume that a good notion of derived categories is given for objects $X$ in 
our $\s$-topos $\T$. Such a good notion includes some version of the base
change property, which is the crucial property to 
pass from the $\s$-operads of spans over an object $X$ to some $\s$-operads of endomorphisms of 
the derived category $\D(X)$. This part uses the full strenght of the Morita theory for dg-categories
as exposed in \cite{to5,to1,tove2}. 

The final section is devoted to a proof of the higher formality conjecture \ref{ci2}. We use the
operations of $\OO$ on $\D(Map(\OO(2),X))$ obtained in the last section in order to
introduce the notion of $\OO$-cohomology of $X$, and prove that two equivalent operads
have equivalent $\OO$-cohomologies. For the particular case of the $E_{k}$-operad (little $k$-disks) 
and the $\mathbb{P}_{k}$-operad ($1-k$-shifted Poisson operad) provides a proof of 
the formality conjecture relating higher Hochschild cohomology and polyvector fields on 
a nice enough derived algebraic stack $X$. The form of the formality theorem is 
an equivalence of dg-Lie algebras, which could be promotted to a more structured
equivalence but this will not be considered in the present work.

\bigskip

\bigskip

\textbf{Acknowledgements:} I am very grateful to A. Blanc, B. Hennion, E. Mann, T. Mignon and M. Robalo
for their participations to the weekly workshop \emph{GW theory and derived algebraic geometry}, 
that has been happening during the first semester 2012-2013 in Montpellier. One of the goal of this
workshop was to clarify how to construct categorified Gromov-Witten invariants as exposed in \cite{seat}. 
The key idea that an operad $\OO$ should act naturally on $\OO(2)$ by cospans, and that this
universal action is responsible for existence of many interesting operations on derived categories, 
has emerged during one of the session. It has been a great surprised to me that the exact
same argument could also be very helpful in deformation quantization and could 
lead to a proof of the higher formality conjecture \ref{cform}. I also thank G. Vezzosi for his
comments on the first versions of this manuscript.

I would like to thank D. Calaque, for many conversations concerning formality, I have learned a lot 
on the subject from him. It seems also that he has another approach to 
the higher formality, based on a generalisation of Tamarkin's proof of Kontsevich's formality. 
This approach however is presently local and not well suited to treat the general case of Artin stacks, even thought 
it is certainly more direct and more explicit. D. Tamarkin and T. Willwacher has also informed me that 
they have made some progress on closely related questions. 

I am thankful to B. Tsygan for conversations on formality and deformation quantization during his visit in Montpellier
in 2011, and for having shared with me some
unpublished notes on equivariant formality. Finally, I would like to thank K. Costello, J. Lurie and T. Preygel
for helpful comments about the Poisson and $E_{k}$-operads. \\

\bigskip

\textbf{Notations and conventions:} We use the expression \emph{$\s$-category} to 
mean \emph{Segal category} (see \cite{si}, see also \cite{tove2}), 
and \emph{strict $\s$-category} to 
mean \emph{simplicial category}. For two $\s$-categories $A$ and $B$ we let 
$Fun^{\s}(A,B)$ be the $\s$-category of $\s$-functors between $A$ and $B$, by which we
mean $\underline{Hom}(A,B')$ for $B'$ a fibrant replacement of $B$ as an object of
the category of Segal precategories, see \cite{si}). 

For a model category $M$ we let $LM$ be its localisation along equivalences, considered
as an $\s$-category (see \cite{tove2}). When $M$ is simplicial we often identify, up to a natural
equivalence, $LM$ with the strict $\s$-category of fibrant and cofibrant objects in $M$. 
When $M=SSet$ is the standard model category of simplicial sets we let $\Top:=LSSet$ be
the $\s$-category of homotopy types. When $M=SeCat$ is the model category of Segal 
precategories we let $\scat:=LSeCat$ be the $\s$-category of $\s$-categories. 

We use $\s$-topos theory from \cite{tove2} (see also \cite{ht}, \cite{seat}). All of our $\s$-topos will
be t-complete by definition (or \emph{hyper-complete} in the terminology of \cite{ht}). For 
an $\s$-topos $\T$ a \emph{stack over $\T$} means a limit preserving $\s$-functor
$\T^{op} \longrightarrow \Top$, where $\Top$ is the $\s$-category of homotopy types. In the same
way, a \emph{categorical stack over $\T$} is a limit preserving $\s$-functor
$\T^{op} \longrightarrow \scat$. 

We denote by $\Omega$ the category of rooted trees of \cite{mw}. For a tree $T\in \Omega$ 
we let $V(T)$ its set of vertices, and for a vertex $v\in V(T)$ we let $T(v)$ be the
star at $v$, that is the maximal subtree of $T$ having $v$ as a unique vertex.
The number of leaves of $T(v)$ is denoted by $n_v$, and is the number of ingoing edges at $v$
(which can be zero). 

\section{Operadic preliminaries}

In this first paragraph we remind notions from operads and $\s$-operads. We start by 
the notion of $\s$-operads according to \cite{cm1,cm2}, defined as certain functors
on a category of rooted trees. We then introduce an auxiliary notion of 
categorical operads, and relate it to pseudo-functors and functors over the
cafegory of trees. 

\subsection{Reminders on (reduced) $\s$-operads}

We will use the expression \emph{$\s$-operad} to mean a \emph{reduced Segal operad} in the sense of 
\cite{cm1,cm2}.
The expression \emph{operad} will rather be restricted to the strict notion of an (uncolored) operad in 
a symmetric mono\"\i dal category. \\

We remind from \cite{mw} 
the category of rooted trees $\Omega$. Its objects are oriented non-planar trees, with 
a unique root and labelled leaves. As in \cite{mw}, we denote by $\eta$ the object of $\Omega$ consisting 
of the unique rooted tree with no vertices. Under the full embedding $\Delta \hookrightarrow \Omega$, the 
object $\eta$ corresponds to the image of $[0] \in \Delta$ (see \cite{mw}). This embedding 
sends a linearly ordered set $[n]=\{0,\dots,n\}$ to the linear tree with 
$n$ vertex and $n+1$ edges.
For a finite set $I$, we denote by $T(I)$ the object of $\Omega$ which is the tree with 
a unique vertex and $card(I)$ leaves labelled by elements of $I$. 
Using the terminology of \cite{mw}, if $I=\{1,\dots,n\}$, then
$T(I)$ is a collar with n leaves and will also be denoted by $T(n)$.  
Note that $T(\emptyset)$ is the
tree with a unique vertex and no leaves and is different from $\eta$.

For a tree $T$ we will generally denote by $V(T)$ its set of vertices. For an vertex
$v \in V(T)$ we denote by $T(v)$ the maximal subtree of $T$ containing the unique
vertex $v$. The tree $T(v)$ is star with $n_v$ leaves, where $n_v$ is the number of ingoing
edges at $v$. \\

The morphisms in $\Omega$ are known to be generated by four kinds of morphisms (see \cite[Lem. 3.1]{mw}).

\begin{enumerate}

\item (automorphisms) Any tree $T\in \Omega$ has a finite automorphism group $aut(T)$.

\item (contraction of an inner edge or inner face) Let $T \in \Omega$ be a tree and $e$ be an inner
edge of $T$ (i.e. an edge bound by two vertices at its end, or equivalently an edge 
which is neither a leaf nor the root). There is a map $T/e \rightarrow T$, called
the contraction along $e$, where $T/e$ consists of the tree where the edge $e$ 
has been contracted to a single vertex. 

\item (Segal or outer face) 
Let $T$ be a tree and $v$ be a vertex for which ingoing edges are all leaves. 
Let $T-v$ be the tree for which the vertex $v$ and all its ingoing leaves have been deleted, and 
$T(v)$ be the sub-tree of $T$ containing $v$ as unique vertex and all its ingoing and outgoing
edges. Then there are morphisms 
$$T-v \longrightarrow T \qquad T(v) \longrightarrow T.$$

\item (identities) For a tree $T$ and any edge $e$ of $T$, let $T+e$ be the tree 
where the edge $e$ has been replaced by the linear tree $[1]$ (two edges and a unique vertex). Then
there is a morphism $T+e \longrightarrow T$. 

\end{enumerate}

The morphisms of type $(2)$ and $(3)$ are called \emph{faces}, and the morphisms 
of type $(4)$ \emph{degeneracies}. According to \cite[Lem. 3.1]{mw} any morphism in $\Omega$ can be
decomposed as an automorphism, a face map and a degeneracy map. \\

The category $\Omega$ can also be considered as a full sub-category of
the category $Op$, of symmetric, reduced and unital operads in the category of sets. 
A full embedding $\Omega \hookrightarrow Op$ is constructed by sending a tree 
$T$ to $\Omega[T]$, the free operad generated by $T$ (see \cite{mw}). For instance, 
the operad $\Omega[T(I)]$ is the free operad generated by single $n$-ary operation, where
$n$ is the cardinality of the set $I$. \\

Let $\T$ be an $\s$-category with finite limits. We consider the $\s$-category of $\s$-functors
$Fun^{\s}(\Omega^{op},\T)$, and we denote by 
$PrOp(\T)$ the full sub-$\s$-category of 
$Fun^{\s}(\Omega^{op},\T)$ consisting of all $\s$-functors $\OO : \Omega^{op} \longrightarrow \T$
sending $\eta$ to a final object in $\T$. The $\s$-category 
$PrOp(\T)$ is by definition the \emph{$\s$-category of (reduced) pre-$\s$-operads in $\T$}. 

An \emph{$\s$-operad in $\T$} is by definition an object $\OO \in PrOp(\T)$ such that 
for all tree $T \in \Omega$ with at least one inner edge, the natural morphism
$$\OO(T) \longrightarrow \OO(Sc(T))$$
is an equivalence in $\T$. Here $Sc(T) \hookrightarrow T$ is the 
\emph{Segal core} of the tree $T$, that is the colimit in $\Omega$ of the
category of subobjects of $T$ either of the form $T(I)$ for some finite set $I$, 
or of the form $\eta$. Because our operads are reduced by 
definition, $\OO(\eta)=*$, another way to state this condition is by saying that for
any tree $T \in \Omega$, the natural morphism
$$\OO(T) \longrightarrow \prod_{v}\OO(T(v)),$$
is an equivalence, where $v$ runs over all vertices of $T$ and $T(v)$ is the
sub-tree of $T$ consisting of the star at $v$ ($v$ as a vertex all its ingoing and outgoing edges). 
The full sub-$\s$-category of $PrOp(\T)$ consisting of $\s$-operads will be denoted by 
$\s-Op(\T)$. 

\begin{df}\label{d1}
For an $\s$-category with finite limits $\T$, \emph{the $\s$-category of $\s$-operads in $\T$} is 
the full sub-$\s$-category of $Fun^{\s}(\Omega^{op},\T)$ consisting of $\s$-functors
sending $\eta$ to $*$ and Segal core to equivalences. It is denoted by 
$\s-Op(\T)$. 

When $\T=\Top$ or $\T=\scat$ is the $\s$-category of spaces or of $\s$-categories we use the following 
notations
$$\s-Op:=\s-Op(\Top) \qquad \s-OpCat:=\s-Op(\scat).$$
\end{df}

In the case the $\s$-category $\T$ is an $\s$-topos, it is of the form $LM$, for a model topos $M$ (e.g. 
$M$ is the model category of simplicial presheaves on a simplicial site, see \cite{hagI}), we have 
a model category $SeOp(M)$ of Segal operads in $M$, whose underlying category is the
category of functors $\Omega^{op} \longrightarrow M$ sending 
$\eta$ to the final object $* \in M$. The model structure on 
the category $SeOp(M)$ is such that fibrant objects are functors $\OO : \Omega^{op} \longrightarrow M$
with levelwise fibrant values and such that for all tree $T$ 
the natural morphism (the Segal face map)
$$\OO(T) \longrightarrow \displaystyle{\prod_{v\in V(T)}}\OO(T(v))$$
is an equivalence in $M$.
We refer to \cite{cm1,cm2} for more details on the model structure on $SeOp(M)$ (the generalisation from
$M=SSet$ to the case where $M$ is the model category of simplicial presheaves on a simplicial
site is straightforward).

We will also consider the model category of strict operads $Op(M)$ in the model category
$M$, considered as a symmetric mono\"\i dal model category for the direct product structure. 
As shown in \cite{cm2}, there is an inclusion functor $Op(M) \longrightarrow SeOp(M)$, 
identifying $Op(M)$ with the full sub-category of $SeOp(M)$ of functors $\OO$
for which all the morphisms $\OO(T) \longrightarrow \OO(Sc(T))$ are isomorphisms. The category
$Op(M)$ also carries a model category structure, for which the equivalences and fibrations are
defined in $SeOp(M)$, so that the inclusion functor $Op(M) \longrightarrow SeOp(M)$ is
a right Quillen functor. By \cite[Thm. 8.15]{cm2} this inclusion is a Quillen equivalence. 

The rectification theorem for diagrams inside $\s$-categories (see \cite[Prop. 4.2.4.4]{ht} and also \cite{tove2})
implies that there is a natural $\s$-functor
$$LSeOp(M) \longrightarrow \s-Op(LM)\simeq \s-Op(\T),$$
which is an equivalence of $\s$-categories. With the previous Quillen equivalence, we do obtain
a string of equivalences of  $\s$-categories
$$LOp(M) \simeq LSeOp(M) \simeq \s-Op(\T).$$
In particular all $\s$-operads in $\T$ can be modeled by strict operads in the model 
category $M$, a fact which will be extremelly useful for us. 

Finally, note that the $\s$-category $\s-Op(\T)$ is also equivalent to the
$\s$-category of stacks on $\T$ with values in $\s$-operads, 
that is $\s$-functors
$$\OO : \T^{op} \longrightarrow L(\s-Op)\simeq \s-Op(\Top),$$
that commutes with limits.
This point of view will allow us to consider variations of the notion of $\s$-operads in $\T$,
the main one being stacks on $\T$ with values in $\s$-operads in $\s$-categories, or in other
words limits preserving $\s$-functors
$$\OO : \T^{op} \longrightarrow \s-Op(\scat).$$

\begin{df}\label{d1'}
Let $\T$ be an $\s$-topos. The \emph{$\s$-category of categorical $\s$-operads
in $\T$} is the $\s$-category of limit preserving $\s$-functors
$$\T^{op} \longrightarrow \s-OpCat.$$
It will be denoted by $\s-OpCat(\T)$.
\end{df}

As $\s$-operads in $\scat$ also is an $\s$-category of $\s$-functors from $\Omega$ to 
$\scat$, the stacks on $\T$ with values in $\s$-operads in $\s$-categories can also be
described as $\s$-functors
$$\OO : \T^{op} \times \Omega^{op} \longrightarrow \scat,$$
which is limits preserving in the first variable, and satisfies the 
Segal core condition in the second variable. \\

Using the notion of fibered and cofibered $\s$-categories (see \cite{tove2}), 
an object $\OO \in
\s-Op(\scat)$, can also be described 
as a fibered $\s$-category $\widetilde{\OO} \longrightarrow \Omega$ satisfying the Segal core
condition. The morphisms 
between $\OO$ and $\OO'$ in $\s-Op(\scat)$ can then be described by the simplicial set
$$Map^{cart}_{\scat/\Omega}(\widetilde{\OO},\widetilde{\OO'}),$$
which is the full sub-simplicial set of the mapping space of the $\s$-category $\scat/\Omega$, 
consisting of $\s$-functors preserving the cartesian morphisms. 

For two $\s$-operads $\OO$ and $\OO'$ in $\scat$, we define the space of 
\emph{very lax morphisms} to be
$$Map_{\s-OpCat}^{vlax}(\OO,\OO'):=Map_{\scat/\Omega}(\widetilde{\OO},\widetilde{\OO'}),$$
where $\widetilde{\OO}$ and $\widetilde{\OO'}$ are the two corresponding fibered $\s$-categories
over $\Omega$. The construction $\OO \mapsto \widetilde{\OO}$ induces a morphism
$$Map_{\s-OpCat}(\OO,\OO') \longrightarrow Map_{\s-OpCat}^{vlax}(\OO,\OO'),$$
and identifies $Map_{\s-OpCat}(\OO,\OO')$ with the union of connected components of 
$Map_{\s-OpCat}^{vlax}(\OO,\OO')$ consisting of $\s$-functors preserving the 
cartesian morphisms. 

We have an intermediate simpicial set
$$Map_{\s-OpCat}(\OO,\OO') \longrightarrow Map_{\s-OpCat}^{lax}(\OO,\OO') \subset 
Map_{\s-OpCat}^{vlax}(\OO,\OO'),$$
where $Map_{\s-OpCat}^{lax}(\OO,\OO')$ consists of $\s$-functors
$\widetilde{\OO} \longrightarrow \widetilde{\OO'}$ over $\Omega$, which preserve 
cartesian morphisms over the Segal (or outer) faces in $\Omega$ (compare with 
the notion of lax mono\"\i dal $\s$-functors of \cite{ha}, Segal faces play the role
of the inert morphisms in $\Omega$). 

\begin{df}\label{d2}
For two objects $\OO,\OO' \in \s-OpCat$, corresponding to two 
fibered $\s$-categories $\widetilde{\OO},\widetilde{\OO'}$ over $\Omega$, the \emph{space of
lax morphisms from $\OO$ to $\OO'$} is defined to be the full sub-simpicial set of 
$Map_{\scat/\Omega}(\widetilde{\OO},\widetilde{\OO'})$ consisting of $\s$-functors preserving the
cartesian morphisms over Segal faces. It is denoted by
$$Map_{\s-OpCat}^{lax}(\OO,\OO').$$ 
The $\s$-categorical operads and lax morphisms form a full sub-$\s$-category 
of $\s$-categories over $\Omega$. It is denoted by 
$\s-OpCat^{lax}$.
\end{df}

Let $a : \OO \longrightarrow \OO'$ be a lax morphism between categorical $\s$-operad, as defined
above and corresponding to an $\s$-functor over $\Omega$
$$\widetilde{a} : \widetilde{\OO} \longrightarrow   \widetilde{\OO'}.$$

The $\s$-functor $\widetilde{a}$ gives in particular $\s$-functors 
$$a_{T} : \OO(T) \longrightarrow \OO'(T),$$
where $\OO(T)$ and $\OO(T')$ respectively are the fibers taken at $T\in \Omega$ of the two
projections  $\widetilde{\OO},\widetilde{\OO'} \rightarrow \Omega$
(the $\s$-functor $a_{T}$ is well defined up to natural equivalences of $\s$-categories).
These $\s$-functors are not functorial in $T$, but for each morphism $u : T \rightarrow T'$ in $\Omega$, 
there is a natural transformation of $\s$-functors
$$h_{u} : u^{*}a_{T'} \Rightarrow a_{T}u^{*} : \OO(T') \rightrightarrows \OO'(T),$$
well defined up to equivalence. The $\s$-functor $a$ over $\Omega$ also controlls the 
higher coherences satisfied by the $\s$-functors $a_{T}$ and the natural 
transformation $h_{u}$. By definition the natural transformations $h_{u}$ are all 
equivalences when $u$ is a Segal face morphism. The mapping space $Map_{\s-OpCat}(\OO,\OO')$
can be identified with the full sub-simplicial set of $Map_{\scat/\Omega}(\widetilde{\OO},
\widetilde{\OO'})$ consisting of all $\s$-functors $a$ such that $h_{u}$ are equivalences for all 
morphisms in $u$. \\
 
\subsection{Operadic categories}

In this paragraph we introduce the notion of a \emph{weak operad in categories}, which will be called
\emph{Op-categories}. For this we first recall that for a group $G$ and a category $\mathcal{C}$, a 
$G$-action on $\mathcal{C}$ consists of the data of a mono\"\i dal (associative)
functor $G \longrightarrow \underline{End}(\mathcal{C})$ which is strictly unital 
(here $G$ is considered as a discrete mono\"\i dal category). 
In short, it consists for all $g \in G$ of an auto-equivalence
$g(-) : x\mapsto g(x)$ of $\mathcal{C}$ (with $e(-)=Id$) together with a family of natural isomorphisms
$$\{\alpha_{g,h} : g(-)\circ h(-) \simeq (gh)(-)\}_{g,h\in G},$$
satisfying the usual cocycle conditions, that for all $x\in \mathcal{C}$, the following diagram
$$\xymatrix{
g(h(k(x))) \ar[rr]^-{g(\alpha_{h,k}(x))} \ar[dd]_-{\alpha_{g,h}(k(x))} & & g(hk(x)) \ar[dd]^-{\alpha_{g,hk}(x)} \\
 & & \\
gh(k(x)) \ar[rr]_-{\alpha_{gh,k}(x)} & & ghk(x)}$$
commutes. A $G$-equivariant functor
between two categories $\mathcal{C}$ and $\mathcal{D}$ with $G$-actions consists of a
functor $F : \mathcal{C} \longrightarrow \mathcal{D}$ together with for every $g\in G$ a 
natural isomorphism $\beta_{g} : F(g(-)) \simeq g(F(-))$, with $\beta_{e}=id$ and compatible
with the isomorphismes $\alpha_{g,h}$ for $\mathcal{C}$ and $\mathcal{D}$. In the sequel we will
never write explicitely the natural isomorphisms $\alpha_{g,h}$ and $\beta_{g}$. 
The categories with $G$-actions and $G$-equivariant functors form a strict $2$-category 
in an obvious way. The expression \emph{natural transformation of $G$-equivariant functors}
will refer to the $2$-morphisms in this $2$-category. \\

We will denote by $FS$ the category of finite sets. For all $I\in FS$ we denote by 
$\Sigma_{I}$ the group of permutations of $I$. For a morphism $u : I \rightarrow J$ in 
$FS$ we denote by $\Sigma_{u}$ the automorphism group of $u$ as an object in the category 
of morphisms in $FS$. We have a canonical isomorphism
$$\Sigma_{u}\simeq \Sigma_{J}\times \prod_{j\in J}\Sigma_{u^-{1}(j)},$$
and canonical morphism of groups  $\Sigma_{u} \longrightarrow \Sigma_{I}$. 
An \emph{Op-category} $\mathcal{C}$ consists of the following data. 

\begin{enumerate}

\item (equivariance) For every finite set $I \in FS$ a category $\mathcal{C}(I)$ with an action 
of the group $\Sigma_{I}=Aut(I)$.
\item (unit) An object $\mathbf{1}_{\mathcal{C}}\in \mathcal{C}(1)$. 
\item (composition) For each map of finite sets $u : I \rightarrow J$,  a 
$\Sigma_{u}$-equivariant functor
$$\mu_{u} : \mathcal{C}(J) \times \prod_{j \in J}\mathcal{C}(u^{-1}(j)) \longrightarrow
\mathcal{C}(I).$$
\item (associativity constraints) For every pair of composable maps of finite sets $\xymatrix{I \ar[r]^-{u} & J \ar[r]^-{v} & K}$
a natural isomorphism of $\Sigma_{vu}$-equivariant functors
$\phi_{u,v} : \mu_{u} \circ \mu_{v} \simeq \mu_{vu} \circ \prod_{k} \mu_{u_{k}}$, 
$$\xymatrix@1{
\mathcal{C}(K) \times \displaystyle{\prod_{k\in K}}\mathcal{C}(v^{-1}(k)) \times  \displaystyle{\prod_{j\in v^{-1}
(k)}}\mathcal{C}(u^{-1}(j))
\ar[rr]^-{\mu_{v} \times id} \ar[d]_-{\prod \mu_{u_{k}}} & & 
\mathcal{C}(J) \times  \displaystyle{\prod_{k \in K}\prod_{j\in v^{-1}(k)}}
\mathcal{C}(u^{-1}(j)) \ar[d]^-{\mu_{u}} \ar@{=>}
[dll]^-{\phi_{u,v}} \\
\mathcal{C}(K) \times  \displaystyle{\prod_{k \in k}\mathcal{C}((vu)^{-1}(k))}
 \ar[rr]_-{\mu_{vu}} & & \mathcal{C}(I),}$$ 
where $u_{k} : (vu)^{-1}(k) \rightarrow v^{-1}(k)$ is the restriction of $u$. 
\item (right unit constraint) For all finite set $I$, a natural isomorphism of $\Sigma_{I}$-equivariant functors
$\alpha_{I} : pr \simeq \mu_{id}\circ id\times (\mathbf{1}_{\mathcal{C}})^{I}$ 
$$\xymatrix{\mathcal{C}(I) \ar[d]_-{id\times (\mathbf{1}_{\mathcal{C}})^{I}} \ar[rd]^-{id} & \\
\mathcal{C}(I) \times \displaystyle{\prod_{i\in I}}\mathcal{C}(\{i\}) \ar[r]_-{\mu_{id}} & \mathcal{C}(I),}$$
from the projection on the first factor $pr$ and $\mu_{id}$. 
\item (left unit constraint) For all finite set $I$, a natural isomorphism of $\Sigma_{I}$-equivariant functors
$\beta_{I} : pr \simeq \mu_{id}\circ id\times (\mathbf{1}_{\mathcal{C}})^{I}$ 
$$\xymatrix{\mathcal{C}(I) \ar[d]_-{\mathbf{1}_{\mathcal{C}}\times id} \ar[rd]^-{id} & \\
\mathcal{C}(*) \times \mathcal{C}(I) \ar[r]_-{\mu_{p}} & \mathcal{C}(I),}$$
where $p : I \rightarrow *$ is the projection to the final object. 
\end{enumerate} 

These data must satisfy some standard cocycle conditions, expressing that the 
natural isomorphisms $\phi_{u,v}$, $\alpha_{I}$ and $\beta_{I}$ are compatible with the composition of maps in 
the category of finite sets $FS$. We leave these compatibilities to the reader. \\

For two Op-categories $\mathcal{C}$ and $\mathcal{D}$, an \emph{Op-functor} 
$f : \mathcal{C} \longrightarrow \mathcal{D}$ consists of the following data.

\begin{enumerate}
\item For every finite set $I$, a $\Sigma_{I}$-equivariant functor
$$f_{I} : \mathcal{C}(I) \longrightarrow \mathcal{D}(I).$$
\item An isomorphism in $\mathcal{D}(*)$ 
$$\gamma_{e} : \mathbf{1}_{\mathcal{C}}\simeq f_{*}(\mathbf{1}_{\mathcal{C}}).$$
\item For every morphism $u : I \rightarrow J$ in $FS$, a natural isomorphism of 
$\Sigma_{u}$-equivariant functors
$$\psi_{u} : \mu_{u}\circ (f_{J}\times \prod_{j\in J}f_{u^{-1}(j)}) \simeq f_{I} \circ \mu_{u} $$
$$\xymatrix@1{
\mathcal{C}(J) \times \displaystyle{\prod_{j\in J}}\mathcal{C}(u^{-1}(j)) 
\ar[rr]^-{\mu_{u}} \ar[d]_-{f_{J}\prod f_{u^{-1}(j)}} & & 
\mathcal{C}(I)  \ar[d]^-{f_{I}} \ar@{<=}
[dll]^-{\psi_{u}} \\
\mathcal{D}(J) \times  \displaystyle{\prod_{j \in j}\mathcal{D}(u^{-1}(j))}
 \ar[rr]_-{\mu_{u}} & & \mathcal{D}(I).}$$ 
\end{enumerate}
These data are required to satisfy some compatibility left to the reader. For two
Op-functors $f$ and $g$, from $\mathcal{C}$ to $\mathcal{D}$ there is an obvious
notion of Op-natural transformations $h : f \Rightarrow g$. It consists for all 
finite set $I$ of a natural transformation $h_{I} : f_{I} \Rightarrow g_{I}$, which 
is compatible with the natural isomorphisms $\gamma_{e}$ and $\psi_{u}$ above. 

The Op-categories, Op-functors and Op-natural transformations form a strict 
$2$-category in an obvious manner. 
It will be denoted by $\Opcat$. In the definition of 
Op-functor the condition that $\gamma_{e}$ and the $\psi_{u}$'s are 
isomorphisms can be relaxed by only requiring that they
exist as non-invertible natural transformations, 
in which case we will talk about \emph{lax Op-functors}. The 
Op-categories, lax Op-functors and Op-natural transformations form another
strict 2-category, denoted by $\Opcat^{lax}$. There is an inclusion of 2-categories
$$\Opcat \subset \Opcat^{lax},$$
identifying $\Opcat$ with the strict sub-2-category consisting of lax Op-functors for which 
$\gamma_{e}$ and the $\psi_{u}$'s are all isomorphisms. \\

We finish this paragraph by describing a rectification functor
$$\Opcat \longrightarrow Fun(\Omega^{op},Cat),$$
from the 1-category of Op-categories to the 1-category of functors from $\Omega^{op}$ 
to categories. For an Op-category $\mathcal{C}$ we define a  pseudo-functor
$$F_{\mathcal{C}} : \Omega^{op} \longrightarrow Cat,$$
by sending a tree $T \in \Omega^{op}$ to the category 
$$F_{\mathcal{C}}(T):=\prod_{v\in V(T)}\mathcal{C}(I_{v}(T)),$$
where $v$ runs over all vertices of $T$ and $I_{v}(T)$ is the valency of $v$ (i.e. the set of
ingoing edges at $v$). Note that if $T$ has no vertex then the formula reads
$F_{\mathcal{C}}(T)=*$.  

The action of $F_{\mathcal{C}}$ on morphisms in $\Omega$ is defined using the
structure functors defining the Op-category $\mathcal{C}$. They can be described
for each of the generating morphism: automorphism, inner face, Segal face and degeneracies
as follows. 

\begin{enumerate}
\item 
For an automorphism
$\tau : T \simeq T$ of a tree $T\in \Omega$, the induced action
$$\tau : \prod_{v}\mathcal{C}(I_{v}(T)) \simeq \prod_{v}\mathcal{C}(I_{v}(T))$$
is defined in a natural way by the $\Sigma_{I_{v}(T)}$-action on 
the category $\mathcal{C}(I_{v}(T))$. 

\item Let $T$ be a tree, $e$ an inner edge of $T$ and 
$T/e \longrightarrow T$ the corresponding morphism in $\Omega$. Let 
$v_{0}$ and $v_{1}$ the two boundary vertices of $e$ in $T$, and 
$x$ the new vertex in $T/e$ created by the contraction. 
Then $V(T/e)=V(T)-\{v_{0},v_{1}\}\coprod\{x\}$. 
The functor
$$c_{e} : \prod_{v \in V(T)}\mathcal{C}(I_{v}(T)) \longrightarrow \prod_{v \in V(T/e)}\mathcal{C}(I_{v}
(T/e)),$$
is defined as follows. For $v \in V(T/e)$ different from $x$, we have
$I_{v}(T/e)=I_{v}(T)$, and $I_{x}(T/e)=I_{v_{0}}(T)\circ_{e} I_{v_{1}}(T)$, the gluing of
the sets $I_{v_{0}}(T)$ and $I_{v_{1}}(T)$ along the element $e \in I_{v_{1}}(T)$. 
The functor $c_{e}$ is defined by a family of functors
$$c_{e,w} : \prod_{v \in V(T)}\mathcal{C}(I_{v}(T)) \longrightarrow \mathcal{C}(I_{w}(T)) \qquad
c_{e,x} : \prod_{v \in V(T)}\mathcal{C}(I_{v}(T)) \longrightarrow \mathcal{C}(I_{v_{0}}(T)\circ 
I_{v_{1}}(T)).$$
For $w\neq x$, the functor $c_{e,w}$ is taken to be the projection. 
The functor
$c_{e,x}$ is defined to be the following composition
$$\xymatrix{
\prod_{v \in V(T)}\mathcal{C}(I_{v}(T)) \ar[r]^-{pr} & 
\mathcal{C}(I_{v_{1}}(T)) \times \mathcal{C}(I_{v_{0}}(T)) \ar[r]^-{\mathbf{1}} & }$$
$$\xymatrix{
\mathcal{C}(I_{v_{1}}(T)) \times \displaystyle{\prod_{j \in I_{v_{1}}(T)}} 
\mathcal{C}(u^{-1}(j)) \ar[r]^-{\mu_{u}} & \mathcal{C}(I_{v_{0}}(T)\circ 
I_{v_{1}}(T)),}$$
where $u : I_{v_{0}}(T)\circ_{e} I_{v_{1}}(T) \longrightarrow I_{v_{1}}(T)$
is the natural projection contracting $I_{v_{0}}(T)$ to $e$. Note that we have
denoted symbolically by
$$\mathbf{1} : \mathcal{C}(I_{v_{1}}(T)) \times \mathcal{C}(I_{v_{0}}(T)) \Longrightarrow
\mathcal{C}(I_{v_{1}}(T)) \times \displaystyle{\prod_{j \in I_{v_{1}}(T)}} 
\mathcal{C}(u^{-1}(j)),$$
the natural functor defined using the units objects $\mathbf{1} \in \mathcal{C}(u^{-1}(j))=\mathcal{C}(\{j\})$
for any $j\neq e$, as well as the identity functors 
of $\mathcal{C}(I_{v_{1}}(T))$ and  $\mathcal{C}(u^{-1}(e))$.

\item Let $T$ be a tree and $T' \longrightarrow T$ a morphism in $\Omega$ corresponding to 
a subtree which is the union of sub-collars in $T$ (i.e. $T'$ is determined by a subset 
$V(T')$ of the set of vertices $V(T)$ of $T$). The functor
$$F_{\mathcal{C}}(T)=\prod_{v\in V(T)}\mathcal{C}(I_{v}(T)) 
\longrightarrow F_{\mathcal{C}}(T')=\prod_{v\in V(T')}\mathcal{C}(I_{v}(T))$$
is simply defined to be the natural projection on the subfactors corresponding
to $V(T') \subset V(T)$. 

\item Let $T$ be a tree and $e$ be an edge. Let $T+e \longrightarrow T$ be the corresponding
morphism in $T$, which doubles the edge $e$ by adding a vertex $x$ to $T$. Let $v_{0}$ and
$v_{1}$ the vertices bouding $e$ in $T$. We have 
$V(T+e)=V(T) \coprod \{x\}$, and 
$$I_{v}(T+e)=I_{v}(T) \; for \; v\neq e \qquad 
I_{e}(T+e)=*.$$
The functor
$$F_{\mathcal{C}}(T)=\prod_{v\in V(T)}\mathcal{C}(I_{v}(T)) \times *  \longrightarrow 
F_{\mathcal{C}}(T+e)=\prod_{v\in V(T)}\mathcal{C}(I_{v}(T)) \times \mathcal{C}(*)$$
is defined to be the product of the identity on the component $\prod_{v\in V(T)}\mathcal{C}(I_{v}(T))$
with the unit object $\mathbf{1} : * \longrightarrow \mathcal{C}(*)$. 
\end{enumerate}

The above rules define $F_{\mathcal{C}}$ on objects and morphisms. The coherences needed to 
promote $F_{\mathcal{C}}$ to a pseudo-functor are themselves naturally defined using the
natural isomorphisms $\phi_{u,v}$, $\alpha_{I}$ and $\beta_{I}$ of the definition 
of an Op-category. We leave this part to the reader. \\

The above construction defines a functor
$$F_{-} : \Opcat \longrightarrow Fun^{ps}(\Omega^{op},Cat).$$
We compose this with the rectification functor of pseudofunctors, in order to 
obtain 
$$\begin{array}{ccc}
\Opcat & \longrightarrow & Fun(\Omega^{op},Cat) \\
 \mathcal{C} & \mapsto  & \mathcal{C}^{str},
\end{array}$$
which turn functorially Op-categories into strict functors from trees to categories. The construction
$\mathcal{C} \mapsto \mathcal{C}^{str}$ will be called the \emph{rectification} or \emph{strictification}
construction. It can shown to induce an equivalence between the 2-category of Op-categories and
the full sub $2$-category of strict functor $\Omega^{op} \longrightarrow Cat$
satisfying the Segal core condition, but this result will not be needed.

\section{Operads of (co)spans}

We let $\T$ be a (t-complete) 
$\s$-topos and $X\in \T$ an object. We will
associate to it a categorical $\s$-operad $\underline{End}^{\sqcup}(X)$ in $\T$. 
This $\s$-operad should be understood as classifying 
families of cospans, also called co-corespondences, between sums of copies of $X$. 

We choose a presentation $\T \simeq LM$ for $M$ a model topos (see \cite{hagI}). We 
can choose $M$ to be a combinatorial, simplicial and proper model category, for which 
the cofibrations are moreover the monomorphisms. By using the injective model
structure on simplicial presheaves on a semi-model site as in \cite{hagI}
we moreover assume that the underlying category of $M$ is a presheaf category.

The model
category $M$ is endowed with the canonical model topology (see \cite{hagI}), 
for which the covering morphisms are the epimorphic family and will therefore
be considered as a model site. The $\s$-category $\T$ is also
equivalent to the $\s$-category of small stacks over the model site $M$. 
We warn the reader that this
model site being not small some care must be taken by either fixing universes as in \cite{hagI}, or
by chosing a small but large enough full sub-category $C$ of $M$ as a generating site (see e.g. the 
proof of \cite[Thm. 4.9.2]{hagI}).

We now choose an object $F\in M$ which is a model for $X \in \T$. 
For all $x\in M$ we define an Op-category
$\mathcal{E}(F)(x)$ as follows.  
For each finite set $I$
the category $\mathcal{E}(F)(x,I)$ possesses as objects the diagrams in $M/x$ 
of the following form
$$\xymatrix{ \displaystyle{\coprod_{I}}F_{|x} \ar[r]& F_{0} & \ar[l] F_{|x},}$$
where $F_{|x}:=F\times x$ as an object of $M/x$.  
Morphisms in the category $\mathcal{E}(F)(x,I)$, from an object $\xymatrix{ \coprod_{I}F_{|x}  \ar[r]&  F_{0}
 & \ar[l] F_{|x}}$
to an object $\xymatrix{ \coprod_{I}F_{|x}  \ar[r]&  F_{1}
 & \ar[l] F_{|x}}$ are defined to be commutative diagrams in $M/x$
$$\xymatrix{
\coprod_{I}F_{|x} \ar[d]_-{id} \ar[r] &  F_{0} \ar[d]&  \ar[l] F_{|x} \ar[d]^-{id} \\
\coprod_{I}F_{|x}  \ar[r] & F_{1} & \ar[l]  F_{|x}. }$$
Note that the category $\mathcal{E}(F)(x,I)$ has a natural model structure induced from
the model structure on $M$: the equivalences and fibrations are
defined on the underlying object $F_{0} \in M$.
This model structure is moreover simplicial, combinatorial and proper. Note that the
cofibrant objects in $\mathcal{E}(F)(x,I)$ are diagrams 
$\xymatrix{ \coprod_{I}F_{|x}  \ar[r]&  F_{0}
 & \ar[l] F_{|x}}$ for which the induced morphism 
$\coprod_{I}F_{|x} \coprod F_{|x} \longrightarrow F_{0}$ is a cofibration of
simplicial presheaves (i.e. a monomorphism).
 
For each finite set $I$, the group $\Sigma_{I}$ acts in an obvious way on the category 
$\mathcal{E}(F)(x,I)$, by acting on the object $\coprod_{I}F_{|x}$ by permutation of the different factors. 
For any map of finite sets $u : I \longrightarrow J$, 
we define a functor
$$\mu_{u} : \mathcal{E}(F)(x,J) \times \prod_{j\in J}\mathcal{E}(F)(x,u^{-1}(j))
\longrightarrow \mathcal{E}(F)(x,I)$$
as follows. This functor associates to an object $\xymatrix{ \coprod_{J}F_{|x}  \ar[r]&  F_{J}
& \ar[l] F_{|x}}$ and a family of objects 
$\{\xymatrix@1{ \coprod_{u^{-1}(j)}F_{|x}  \ar[r]&  F_{j} & \ar[l] F_{|x}} \}_{j \in J}$
the object
$$\coprod_{j\in J} \coprod_{u^{-1}(j)}F_{|x} \longrightarrow \mu_{u}(F) \longleftarrow  
F_{|x},$$
defined by the push-out square in $M$
on $C/x$
$$\xymatrix@1{
\displaystyle{\coprod_{I}F_{|x}=\coprod_{j \in J} \coprod_{u^{-1}(j)}F_{|x}} \ar[r] & 
\displaystyle{\coprod_{j\in J}F_{j}}\ar[d] & \displaystyle{\coprod_{J}F_{|x}} \ar[l] \ar[d] \\
& \mu_{u}(F) & F_{J} \ar[l] & F_{|x}. \ar[l]}$$

The unit for the Op-category $\mathcal{E}(F)(x)$ is
the diagram of identities $\xymatrix{F_{|x} & F_{|x} \ar[r] \ar[l] &F_{|x}}$, considered
as an object in $\mathcal{E}(F)(x,*)$. 
We leave to the reader the definition of the natural isomorphisms $\phi_{u,v}$, 
$\beta_{I}$ and $\alpha_{I}$ making the family 
of categories $\mathcal{E}(F)(x,I)$ into an Op-category, induced from the
natural associativity constraints of coproducts in $M$.
The operation functors $\mu_{u}$ are moreover compatible with the model structures
on the categories $\mathcal{E}(F)(x,I)$, in the sense that they are simplicial left Quillen functors, and thus
preserve cofibrations, trivial cofibrations and equivalences between cofibrant objects. \\

When $I$ is fixed, and $u : y \rightarrow x$ is a morphism in $M$, we have a base change functor
$$u^{*} : \mathcal{E}(F)(x,I) \longrightarrow \mathcal{E}(F)(y,I),$$
induced by the base change $M/x \longrightarrow M/y$ sending $z \rightarrow x$ to $z\times_{x}y\rightarrow 
y$. The base change
functors are right Quillen functors, so preserve fibrations, trivial fibrations and equivalences
between fibrant objects.
It is also naturally a simplicial functor. We have
a natural isomorphism $(uv)^{*} \simeq v^*u^*$ for two 
morphisms $\xymatrix{z \ar[r]^-{v} & y \ar[r]^-{u} & x}$
in $M$, making $x \mapsto \mathcal{E}(F)(x,I)$ into a pseudo-functor 
from $M^{op}$ to model categories. Moreover, the functors $u^{*}$ naturally 
commutes with push-outs in $M$ (as $M$ is a presheaf category) and thus the association
$x \mapsto \mathcal{E}(F)(x,-)$ defines a pseudo-functor $M^{op} \longrightarrow \Opcat$ from $M^{op}$ 
to the category of Op-categories. Using the rectification 
$\Opcat \longrightarrow Fun^{ps}(\Omega^{op},Cat)$
we will consider this as a pseudo-functor in two variables
$$\mathcal{E}(F)(-,-) : M^{op} \times \Omega^{op} \longrightarrow Cat.$$
This is compatible with the simplicial enrichements of the categories 
$\mathcal{E}(F)(x,I)$, and thus can be promoted to a pseudo-functor
from $M^{op} \times \Omega^{op}$ to the 2-category of strict $\s$-categories. 
For a morphism of trees $f : T \rightarrow T'$ in $\Omega$ we will denote by 
$f^{*} : \mathcal{E}(F)(-,T') \longrightarrow \mathcal{E}(F)(-,T)$ the corresponding natural transformation
of pseudo-functors on $M^{op}$, which is now a simplicial functor. 

We now define a strict $\s$-category $\widetilde{\mathcal{E}(F)}$, over $C\times \Omega$ as follows
by applying the Grothendieck's construction  to the pseudo-functor $\mathcal{E}(F)(-,-)$.
Its objects are triples $(x,T,c)$, where $(x,T) \in M\times \Omega$, and $c$ is an object 
in $\mathcal{E}(x,T)$. For two objects $(x,T,c)$ and $(y,T',c')$, we define a simplicial set of maps
$$\widetilde{\mathcal{E}(F)}((x,T,c),(y,T',c')) := \coprod_{(u,f) \in M(x,y) \times \Omega(T,T')}
\underline{Hom}(u^{*}f^*(c'),c),$$
where $\underline{Hom}(u^{*}f^{*}(c'),c)$ is the simplicial hom object of the simplicial 
category $\mathcal{E}(F)(x,T)$. Morphisms are composed in the obvious way making 
$\widetilde{\mathcal{E}(F)}$ into a simplicial category (or a strict $\s$-category). It comes
equiped with a natural projection
$$\pi : \widetilde{\mathcal{E}(F)} \longrightarrow M\times \Omega,$$
making it into a fibered $\s$-category. The corresponding $\s$-functor simply sends 
$(x,T)$ to the
simplicial category $\mathcal{E}(F)(x,T)$, considered as a strict $\s$-category, and 
a morphism $(u,f) \in M(x,y) \times \Omega(T,T')$ to the strict $\s$-functor $u^{*}f^{*}\simeq f^{*}u^{*}$.

We let $\widetilde{\mathcal{E}(F)}^{cf}$ be the full sub-$\s$-category of $\widetilde{\mathcal{E}(F)}$ 
consisting
of objects $(x,T,c)$ with $c$ a cofibrant and fibrant object in $\mathcal{E}(F)(x,T)$.

\begin{lem}\label{l1}
The $\s$-category $\widetilde{\mathcal{E}(F)}^{cf}$ is fibered over $M \times \Omega$. 
\end{lem}

\textit{Proof:} This is a general fact. Let $I$ be a small category and 
$\mathcal{F} : I^{op} \longrightarrow  \mathbb{S}-Cat$
a pseudo-functor from $I^{op}$ to the 2-category of 
simplicial categories (or strict $\s$-categories). Assume that for all $i \in I$ the category $\mathcal{F}(i)$ is endowed with 
a simplicial model structure. Assume moreover that for all morphism $u : i \rightarrow j$ in $I$, 
the base change functor $u^{*} : \mathcal{F}(j) \longrightarrow \mathcal{F}(i)$ 
preserves cofibrations as well as equivalences between 
fibrant and cofibrant objects. We let $\widetilde{\mathcal{F}}$ be the (strict)  $\s$-category of pairs 
$(i,x)$, with $x \in \mathcal{F}(i)$. 
The simpicial set of maps $(i,x) \rightarrow (j,y)$ in $\widetilde{\mathcal{F}}$ is defined to be
$$\widetilde{\mathcal{F}}((i,x),(j,y))=\coprod_{u \in I(i,j)}\underline{Hom}(u^{*}(y),x).$$
There is an obvious projection $\widetilde{\mathcal{F}} \longrightarrow I$. 
Let $\widetilde{\mathcal{F}}^{cf} \subset 
\widetilde{\mathcal{F}}$
be the full sub-$\s$-category of pairs $(i,x)$ with $x$ a cofibrant and fibrant object
in $\mathcal{F}(i)$. Then the lemma is a special case
of the following more general statement.

\begin{slem}\label{sl1}
With the above notations the $\s$-category 
$\widetilde{\mathcal{F}}^{cf}$ is fibered over $I$. 
\end{slem}

\textit{Proof of the sub-lemma:}
Let $u : i\rightarrow j$ be a morphism in $I$, and $y \in \mathcal{F}(j)$ be a fibrant and cofibrant 
object. Let 
$\alpha : u^{*}(y) \longrightarrow x$ be a fibrant replacement in $\mathcal{F}(i)$ 
(so $x$ is fibrant and cofibrant
by hypothesis on $u^{*}$). 
The morphism $\alpha$ defines a morphism in $(x,i) \rightarrow (y,j)$ in $\widetilde{\mathcal{F}}^{cf}$
over $u$, and which is easily checked to be a cartesian lift of $u$ to $\widetilde{\mathcal{F}}^{cf}$. 
\hfill $\Box$ \\

This finishes the proof of the lemma \ref{l1}. \hfill $\Box$ \\

By lemma \ref{l1} and \cite[Prop. 1.4]{tove2} we can turn the fibered $\s$-category 
$\widetilde{\mathcal{E}(F)}^{cf}$ into an $\s$-
functor
$$\underline{End}^{\sqcup}(F) : M^{op}\times \Omega^{op} \longrightarrow \scat.$$
The value at $(x,T)$ is naturally equivalent to the strict $\s$-category $\mathcal{E}(F)(x,T)^{cf}$, 
of fibrant and cofibrant objects in $\mathcal{E}(F)(x,T)$. This is also equivalent to the localisation 
$L\mathcal{E}(F)(x,T)$ of the category $\mathcal{E}(F)(x,T)$ along its subcategory of equivalences. 
For a morphism $(u,f) : (x,T) \rightarrow (y,T')$ in $C\times \Omega$
the base change along $(u,f)$ is then naturally equivalent to the derived $\s$-functor
of the pull-back functor  $f^{*}u^{*} : \mathcal{E}(F)(y,T') \longrightarrow \mathcal{E}(F)(x,T)$
$$Lf^{*}Ru^{*} :\xymatrix{ L\mathcal{E}(F)(y,T') & \ar[l]_-{\sim}  L\mathcal{E}(F)(y,T')^{cf} 
\ar[r]^-{f^{*}u^{*}}
& L\mathcal{E}(F)(x,T)^{c} \ar[r]^-{\sim} & L\mathcal{E}(F)(x,T).}$$

The $\s$-functor $\underline{End}^{\sqcup}(F)$ also
possesses the following two exactness properties.

\begin{enumerate}

\item For a fixed object $M\in M$, the $\s$-functor 
$$\underline{End}^{\sqcup}(F)(x,-) : \Omega^{op} \longrightarrow \scat$$
satisfies the Segal core conditions and thus defines an $\s$-categorical operad.

\item For a fixed tree $T \in \Omega$, the $\s$-functor
$$\underline{End}^{\sqcup}(F)(-,T) : M^{op} \longrightarrow \scat$$
sends homotopy colimits in $M$ to homotopy limits, and thus defines a
stack of $\s$-categories over the model site $M$, or equivalently 
an object in $\s-Cat(\T)$. 

\end{enumerate}

Property $(1)$ holds by construction: for a tree $T$, we have a natural equivalence
$$\underline{End}^{\sqcup}(F)(x,T)\simeq \prod_{v\in V(T)}\underline{End}^{\sqcup}(F)(x,T(n_{v})).$$
The property $(2)$ above follows directly from the fact that the prestack of stacks, 
sending $x\in \T$ to $\T/x$, is itself a categorical stack (see \cite[A.2]{to1}).

The two above properties imply that $\underline{End}^{\sqcup}(F)$ defines a categorical $\s$-operad over the 
$\s$-topos $\T$. It is easy to see that as an object in $\s-OpCat(\T)$ it does not depend, up to
equivalences, on the choice of the model topos $M$ and the model $F \in M$. We will therefore denote
it by $\underline{End}^{\sqcup}(X) \in \s-OpCat(\T)$. 

\begin{df}\label{d3}
Let $\T$ be an $\s$-topos and $X\in \T$ an object.
The \emph{$\s$-operad of co-spans over $X$ is} $\underline{End}^{\sqcup}(X) \in \s-OpCat(\T)$.
\end{df}

The $\s$-operad $\underline{End}^{\sqcup}(X)$ is functorial in $X$ in the following sense.
Let $X \rightarrow Y$ be a morphism in $\T$, and choose models $F \rightarrow F'$ in a model
topos $M$ such that $LM\simeq \T$ as above. 

Recall the two pseudo-functors
$$\mathcal{E}(F)(-,I), \mathcal{E}(F')(-,I) : M^{op} \longrightarrow \Opcat.$$
For all finite set $I$, we define a natural transformations 
$$f_{I} : \mathcal{E}(F)(-,I) \longrightarrow \mathcal{E}(F')(-,I),$$
between pseudo-functors as follows.
For $x\in M$, $f_{I}(x)$ sends 
an object $\xymatrix{ \displaystyle{\coprod_{I}}F_{|x} \ar[r]& F_{0} & \ar[l] F_{|x},}$
in $\mathcal{E}(F)(x,I)$
to $\xymatrix{ \displaystyle{\coprod_{I}}F'_{|x} \ar[r]& F'_{0} & \ar[l] F'_{|x}} \in \mathcal{E}(F')(x,I)$,
where $F_{0}'$ is defined to by the push-out in $M$
$$F'_{0}:=\left(\displaystyle{\coprod_{I}}F'_{|x} \coprod F'_{|x} \right) 
\coprod_{\left(\displaystyle{\coprod_{I}}F_{|x} \coprod F_{|x}\right)} F_{0}.$$
Using the functorial associativity and unit isomorphisms for push-outs, the 
family of functors $f_{I}$ defines an natural transformation of pseudo-functors of Op-categories 
on $M^{op}$
$$f : \mathcal{E}(F) \longrightarrow \mathcal{E}(F').$$
The natural transformation is also compatible with the simplicial structures on 
the model categories $\mathcal{E}(F)(x,I)$ and $\mathcal{E}(F')(x,I)$, and thus 
defines a fibered strict $\s$-category
$$\widetilde{\mathcal{E}(f)} \longrightarrow \Delta^{1} \times M \times \Omega,$$
where $\Delta^{1}$ is the category with two objects $0$ and $1$, and unique non trivial morphism
$0 \rightarrow 1$. The general proof of the lemma \ref{l1} shows that the full sub-$\s$-category
$\widetilde{\mathcal{E}(f)}^{cf}$ of cofibrant and fibrant objects remains fibered over 
$\Delta^{1} \times M \times \Omega$. We thus get a well defined  
$\s$-functor of $\s$-categories over $M\times \Omega$
$$f : \widetilde{\mathcal{E}(F)}^{cf} \longrightarrow \widetilde{\mathcal{E}(F')}^{cf}.$$
This $\s$-functor is easily checked to be cartesian and thus
defines a morphism
of categorical $\s$-operads in $\T$
$$f : \underline{End}^{\sqcup}(X) \longrightarrow \underline{End}^{\sqcup}(X').$$
With a bit of care this construction can be naturally lifted to an $\s$-functor
$$\underline{End}^{\sqcup}(-) : \T \longrightarrow \s-OpCat(\T).$$ 

The categorical $\s$-operads $\underline{End}^{\sqcup}(X)$ have a dual version, using 
spans instead of cospans. As above, let $\T$ be of the form $LM$, for a model topos $M$.
Let $F \in M$ be a model for $X \in \T$, which we choose to 
be fibrant.

For $x\in M$, 
We define an Op-category
$\mathcal{E}^{o}(F)(x)$ as follows. For a finite set $I$, the objects of the category 
$\mathcal{E}^{o}(F)(x,I)$ consists of diagrams in $M/x$
$$\xymatrix{ \displaystyle{\prod_{I}}F_{|x} & \ar[l] F_{0}  \ar[r] & F_{|x}.}$$
The morphisms in $\mathcal{E}^{o}(F)(I)$ from $\xymatrix{ \displaystyle{\prod_{I}}F_{|x} & \ar[l] F_{0}  \ar[r] & F_{|x}}$ to $\xymatrix{ \displaystyle{\prod_{I}}F_{|x} & \ar[l] F_{1}  \ar[r] & F_{|x}}$,
are the commutative diagrams (note the oppositive convention compare to 
morphisms in $\mathcal{E}^{o}(F)(x,I)$)
$$\xymatrix{
\prod_{I}F_{|x} \ar[d]_-{id}  & \ar[l]  F_{0}  \ar[r] &   F_{|x} \ar[d]^-{id} \\
\prod_{I}F_{|x}  & \ar[l] F_{1} \ar[u] \ar[r] &  F_{|x}.}$$
The groups 
$\Sigma_{I}$ naturally act of the categories $\mathcal{E}^{o}(F)(I)$ by permuting the
components of $\prod_{I}F_{|x}$. The family of categories $\mathcal{E}^{o}(F)(x,I)$ can be 
organized in an Op-category in a similar fashion as for $\mathcal{E}(F)(x,I)$, but for 
which for any map of finite sets $u : I\rightarrow J$ the
$\Sigma_{u}$-equivariant functors
$$\mu_{u} : \mathcal{E}^{o}(F)(x,J) \times \prod_{j\in J}\mathcal{E}^{o}(F)(x,u^{-1}(j))
\longrightarrow \mathcal{E}^{o}(F)(x,I)$$
are defined using pull-backs instead of push-outs in $M$. We get this way 
a lax functor
$$\mathcal{E}^{o}(F)(-,I) : M^{op} \longrightarrow \Opcat,$$
where for a morphism $u : y \rightarrow x$ in $M$ we use the pull-back functor
$u^{*} : M/x \longrightarrow M/y$ in order to define an Op-functor
$$u^{*} : \mathcal{E}^{o}(F)(x,I) \longrightarrow \mathcal{E}^{o}(F)(y,I).$$
This defines a pseudo-functor of Op-categories 
$\mathcal{E}^{o}(F) : M^{op} \longrightarrow \Opcat$. Using the rectification
of Op-categories of \S 1.2, we can
consider $\mathcal{E}^{o}(F)$ as a pseudo-functor
$$\mathcal{E}^{o}(F) : M^{op}\times \Omega^{op} \longrightarrow Cat.$$
Each category $\mathcal{E}^{o}(F)(x,T)$ has moreover a model category structure induced from the
model structure on $M^{op}$ (for which fibrations are monomorphisms, the cofibrations are the
fibrations in $M$, and equivalences
are defined on the underlying object in $M$: the funny 
switch between fibrations and cofibrations is induced by the opposite convention
on morphisms in $\mathcal{E}^{o}(F)(x,T)$). These model structures are simplicial
and proper (and are opposite to combinatorial model categories), and these structures are 
compatible with the pseudo-functor $\mathcal{E}^{o}(F)$ above.

We can from a strict $\s$-category over $M\times \Omega$
$$\widetilde{\mathcal{E}^{o}(F)} \longrightarrow M\times \Omega,$$
whose objects are triples $(x,T,c)$, where $(x,T) \in M\times \Omega$, and $c$ is an object 
in $\mathcal{E}^{o}(x,T)$. For two objects $(x,T,c)$ and $(y,T',c')$, we define a simplicial set of maps
$$\widetilde{\mathcal{E}(F)}((x,T,c),(y,T',c')) := \coprod_{(u,f) \in C(x,y) \times \Omega(T,T')}
\underline{Hom}(u^{*}f^{*}(c'),c),$$
where $\underline{Hom}(u^{*}f^{*}(c'),c)$ is the simplicial hom object of the simplicial 
category $\mathcal{E}^{o}(x,T)$.  

As for the $\s$-operad of cospans we consider the full sub-$\s$-category
$$\widetilde{\mathcal{E}^{o}(F)}^{cf} \subset \widetilde{\mathcal{E}^{o}(F)}$$
consisting of cofibrant and fibrant objects.
By the general proof of lemma \ref{l1}  we observe that the induced projection
$$\widetilde{\mathcal{E}^{o}(F)}^{cf} \longrightarrow M \times \Omega$$
remains fibered and thus defines an $\s$-functor
$$\underline{End}^{\times}(F) : M^{op}\times \Omega^{op} \longrightarrow \scat.$$
Its value at $(x,T)$ is naturally equivalent to $L\mathcal{E}^{o}(F)(x,T)$, the localisation 
of $\mathcal{E}^{o}(F)(x,T)$ along its equivalences, and the base change $\s$-functors are 
given the derived $\s$-functors of $u^{*}$ and $f^{*}$. As for cospans this $\s$-functor
satisfies the following two properties.

\begin{enumerate}

\item For a fixed object $x\in M$, the $\s$-functor 
$$\underline{End}^{\times}(F)(x,-) : \Omega^{op} \longrightarrow \scat$$
satisfies the Segal core conditions and thus defines an $\s$-categorical operad.

\item For a fixed tree $T \in \Omega$, the $\s$-functor
$$\underline{End}^{\times}(F)(-,T) : M^{op} \longrightarrow \scat$$
sends homotopy colimits in $M$ to homotopy limits and thus defines
a categorical stack over $\T$.

\end{enumerate}

The $\s$-functor $\underline{End}^{\times}(F)$ therefore defines a categorical $\s$-operad
over $\T$. Again this does not depend, up to equivalence, on the choice of the model topos
$M$ nor of the model $F$, and will be denoted by $\underline{End}^{\times}(X)$.

\begin{df}\label{d3'}
Let $\T$ be an $\s$-topos and $X\in \T$ an object.
The \emph{$\s$-operad of spans over $X$ is} $\underline{End}^{\times}(X) \in \s-OpCat(\T)$.
\end{df}

As for the operads of cospans, the construction $X \mapsto \underline{End}^{\times}(X)$ 
can be promoted to an $\s$-functor
$$\T^{op} \longrightarrow \s-OpCat(\T).$$
This $\s$-functor sends a morphism $u : X \rightarrow X'$ to the morphism
$u^{*} : \underline{End}^{\times}(X') \longrightarrow \underline{End}^{\times}(X)$ which 
is obtained by pulling-back the spans over $X$ along the morphism $u$. \\

Assume now that $X,Z \in \T$ are two objects. We let $Z^{X} \in \T$ be the  internal hom object
of maps from $X$ to $Z$. We are going to construct a natural morphism of $\s$-categorical operads over $\T$
$$\underline{End}^{\sqcup}(X) \longrightarrow \underline{End}^{\times}(Z^{X}),$$
well defined up to equivalence. For this we choose again a presentation of $\T$ by a model topos $M$, 
and we choose two models $F$ and $G$ in $M$ of 
$X$ and $Z$.
We suppose that $G$ is a fibrant object. The model category $M$ has internal Hom's, which we denote by 
exponentiation, and is more precisely an internal model category (see \cite{to3}).
We can therefore consider the simplicial functor $(G)^{-} : M\longrightarrow M^{op}$, 
sending $x$ to $G^{x}$. For $x\in M$, we 
also have the
restriction $G_{|x}:=G\times x \in M/x$, 
which remains a fibrant object in $M/x$, 
and we have the simplicial endofunctor over $x$
$$(G_{|x})^{-} : M/x\longrightarrow (M/x)^{op}$$ 
where exponentiations are now
taken in $M/x$.

For $(x,T) \in M\times \Omega$, the functor $(G_{|x})^{-}$ induces a simplicial functor
$$(G_{|x})^{-} : \mathcal{E}(F)(x,T) \longrightarrow \mathcal{E}^{o}(G^{F})(x,T),$$
simply sending an object
$$\xymatrix{ \displaystyle{\coprod_{I}}F_{|x} \ar[r]& F_{0} & \ar[l] F_{|x},}$$ to 
$$\xymatrix{ \displaystyle{\prod_{I}}(G_{|x})^{F_{|x}} & \ar[l] (G_{|x})^{F_{0}} \ar[r] & 
(G_{|x})^{F_{|x}}.}$$
Note here that we use the canonical isomorphisms $(G^{F})_{|x} \simeq G_{|x}^{F_{|x}}$.
This functor commutes, up to natural isomorphisms, with the pull-backs functors along morphisms in $M\times 
\Omega$, and thus
induces a well defined $\s$-functor between strict fibered $\s$-categories over $M\times \Omega$
$$\widetilde{\mathcal{E}(F)} \longrightarrow \widetilde{\mathcal{E}^{o}(G^{F})}.$$
Finally, as $G$ has been chosen to be fibrant it is easy to see that the above $\s$-functor of 
fibered $\s$-categories restricts to an $\s$-functor on fibrant and cofibrant objects
$$(G_{|x})^{-} : \widetilde{\mathcal{E}(F)}^{cf} \longrightarrow \widetilde{\mathcal{E}^{o}(G^{F})}^{cf}.$$
It is easy to check that this morphism preserves cartesian morphisms, and thus defines a morphism of 
categorical $\s$-operads
over $\T$
$$(G)^{-} : \underline{End}^{\sqcup}(F) \longrightarrow \underline{End}^{\times}(G^{F}).$$
This construction is independent of the choices of $M$, $F$ and $G$, and will simply 
be denoted by 
$$ \underline{End}^{\sqcup}(X) \longrightarrow \underline{End}^{\times}(Z^{X}).$$
It is well defined, up to equivalences, as a morphism in the $\s$-category $\s-CatOp(\T)$. 

\begin{df}\label{d4}
Let $\T$ be an $\s$-topos and $X,Z\in \T$ two objects.
The \emph{duality map with respect to $Z$} is the morphism of categorical $\s$-operads in $\T$
$$\mathbb{D}_{Z} : \underline{End}^{\sqcup}(X) \longrightarrow \underline{End}^{\times}(Z^{X})$$
constructed above.
\end{df}

\section{Adjoint and co-adjoint representations of $\s$-operads}

The purpose of this section is to prove our first main theorem. It can be summarized
as follows. Recall that for an $\s$-operad $\OO$ we denote by $\OO(n)$ the 
object $\OO(T([n]))$, where $T([n])$ is the tree with a unique vertex and $n$ leaves. 

\begin{thm}\label{t1}
Let $\OO$ be an $\s$-operad in an $\s$-topos $\T$ (assumed t-complete). We suppose that 
we have $\OO(0)\simeq \OO(1)\simeq *$. Then there exists a lax morphism 
of categorical $\s$-operads over $\T$
$$ad : \OO \longrightarrow \underline{End}^{\sqcup}(\OO(2)).$$
The evaluation of this morphism at the tree $T([n])$ can be described by the following 
cospans in $\T/\OO(n)$
$$\xymatrix{
\displaystyle{\coprod_{n}}\OO(n) \times \OO(2) \ar[r] \ar[dr] & \OO(n+1) \ar[d] & \ar[l] \OO(2) \times \OO(n) 
\ar[ld] \\
 & \OO(n), & }$$
given by the natural composition morphisms $\OO(n) \times \OO(2) \rightarrow \OO(n+1)$, as well
as the natural projection $\OO(n+1) \rightarrow \OO(n)$ obtained by composing with the canonical
$0$-ary operation.
\end{thm}

The construction of the morphism $ad$ will be described in several steps above. 
Before getting to the construction and the proof we add a few words concerning the statement itself. 
The $\s$-operad $\OO$ is by assumption
an $\s$-operad in $\T$, or equivalently a limit preserving $\s$-functor
$$\OO : \T^{op} \longrightarrow \s-Op.$$
We consider it as a categorical $\s$-operad over $\T$ by considering the full embedding
$\Pi_{\s} : \Top \hookrightarrow \scat$, which sends a simplicial set to its
$\s$-fundamental groupoid. This $\s$-functor induces another full embedding
when passing to $\s$-operads
$$\Pi_{\s} : \s-Op=\s-Op(\Top) \hookrightarrow \s-OpCat=\s-Op(\scat),$$
which allows us to consider $\OO$ also a categorical $\s$-operad over $\T$. We should have written
$\Pi_{\s}(\OO)$ in the statement \ref{t1}, but for simplicity 
we will omit this notation and will also consider that the construction 
$\Pi_{\s}$ has been applied if necessary. A typical example will be for objects of $T$ themselves: they represent limit 
preserving $\s$-functors $Map(-,x) : \T^{op} \longrightarrow \Top$ and thus
will also be considered as limit preserving $\s$-functors $\T^{op} \longrightarrow \scat$
(i.e. as categorical stack over $\T$). The categorical stack represented by $x$ will still be denoted by $x$.

Another comment concernes the description of the evaluation at a tree of the form $T([n])$. 
We use here implicitely the following fact: for any two objects $x,y \in \T$, the $\s$-category
of morphisms
$Map(x,\underline{End}^{\sqcup}(y)(T(n)))$ is naturally equivalent to the $\s$-category of diagrams
in $\T$ of the form
$$\xymatrix{
\displaystyle{\coprod_{n}}x \times y \ar[r] \ar[dr] & F \ar[d] & \ar[l] y \times x
\ar[ld] \\
 & x, & }$$
or in other words to the double comma $\s$-category  $z/(\T/x)$, where
$z=\displaystyle{\coprod_{n+1}}x \times y$. This explains the description of the evaluation 
of the lax morphism $ad$ at $T(n)$: it is a morphism of categorical stacks over $\T$
$$\OO(n) \longrightarrow \underline{End}^{\sqcup}(\OO(2))(T(n)),$$
and thus can be described as a cospan over $\OO(n)$, between $\coprod_{n}\OO(2) \times \OO(n)$ 
and $\OO(2) \times \OO(n)$.

The morphisms in the diagram
$$\xymatrix{
\displaystyle{\coprod_{n}}\OO(n) \times \OO(2) \ar[r] \ar[dr] & \OO(n+1) \ar[d] & \ar[l] \OO(2) \times 
\OO(n) 
\ar[ld] \\
 & \OO(n), & }$$
can also be described in more details as follows. Let $<k>=\{1,\dots,k\}$. We consider the natural 
inclusion $<n> \rightarrow <n+1>$, which induces
a composition $\s$-functor 
$$\OO(n+1)\simeq \OO(n+1) \times \OO(*)^{n} \times \OO(0)
 \longrightarrow \OO(n).$$
This is the map which \emph{forget the last leaf} in a tree with leaves labelled by $<n+1>$. Note that
this map exists only because we have imposed the condition $\OO(1)=\OO(0)=*$. 

The various morphisms 
$\coprod_{n}\OO(n) \times \OO(2) \longrightarrow \OO(n+1)$ of the left hand side are 
constructed similarly. Let $i\in <n>$ and $<n+1> \longrightarrow  <n>$ be the natural morphism
sending $j\leq i$ to $j$ and $j>i$ to $j-1$ (so $i$ is hit two times). It induces a composition 
$\s$-functor
$$\OO(n)\times \OO(2) = \OO(n)\times \prod_{1\leq j<i}\OO(*) \times \OO(2)\times \prod_{i<j\leq n}\OO(*)  
\longrightarrow \OO(n+1).$$
This morphism has to be precomposed with the transposition $\tau_{i+1,n+1}$
of $\OO(n+1)$ to obtain the $\s$-functor
$$\OO(n) \times \OO(2) \longrightarrow \OO(n+1)$$
which is the i-th component of the left hand side morphism of the diagram above. The right hand 
side is defined similarly by using the morphism $<n+1> \longrightarrow <2>$ sending $n+1$ to $2$ and 
everything else to $1$. \\

Using the language of trees, the projection $\OO(n+1) \rightarrow \OO(n)$ 
corresponds to the following diagram of trees (or rather presheaves over $\Omega$)
$$\xymatrix{
T(n) \ar[r]^-{a} & T & \ar[l]_-{b} T(n+1) \coprod T(0),}$$
where $T$ is the tree with two vertices obtained by gluing $T(0)$ on the last leaf of $T(n+1)$, 
the morphism $a$ is the contraction of the unique inner edge of $T$ and the map
$b$ is the Segal core decomposition of $T$ at its two vertices. In the same way, the
morphism
$$\displaystyle{\coprod_{n}}\OO(n) \times \OO(2) \longrightarrow \OO(n+1),$$
is given by $n$ morphisms of trees
$$\xymatrix{T(n+1) \ar[r]^-{a_{i}} & T & \ar[l]_-{b_{i}} T(n) \coprod T(2),}$$
where $T$ is the tree obtained by gluing $T([2])$ on the i-th leaf of $T(n)$
(with a relabelling of the leaves such that the second leaf of $T(2)$
becomes the $n+1$-th leaf of $T$), 
$a_{i}$ is the contraction of the unique inner edge of $T$, and $b_{i}$ is the Segal core
decomposition of $T$. Finally, the morphism of the right hand side
$$\OO(n+1) \longleftarrow \OO(n) \times \OO(2)$$
is given by the morphism of trees
$$\xymatrix{T(n+1) \ar[r]^-{a} & T & \ar[l]_-{b} T(n) \coprod T(2),}$$
where $T$ is the tree obtained by gluing $T(n)$ on the first leaf of 
$T(2)$, $a$ is the contraction of the unique inner edge of $T$ and $b$ is the Segal core
decomposition. \\

We now start the proof of the theorem \ref{t1}. \\

\textbf{Step 1:} We start by some preparations as well as some notations. We present the $\s$-topos $\T\simeq
LM$ by a model topos $M$. We assume as before that $M$ is a combinatorial, proper, simplicial
internal model category, and a presheaf category. We use the injective model structures on 
presheaves and thus assume that cofibrations in $M$ are the monomorphisms.

By the rectification theorem of \cite{cm2}, the 
$\s$-operad $\OO$ can be chosen to represented by a strict operad $\OO$ in the category $M$ and 
we can arrange things such that $\OO(0)=\OO(1)=*$. 
We let $K:=\OO(2) \in M$, and we consider $\widetilde{\mathcal{E}(K)}$
the strict $\s$-category considered during the definition of the operad
of cospans in $\T$. 
Recall that the objects of $\widetilde{\mathcal{E}(K)}$ are the triples $(x,T,c)$, where
$x \in M$, $T \in \Omega$, and $c$ is an object in 
$\mathcal{E}(K)(x,T)$. The sub-$\s$-category $\widetilde{\mathcal{E}(K)}^{cf}$ consists
of triples $(x,T,c)$ with $c$ cofibrant and fibrant in $\mathcal{E}(K)(x,T)$. 
We thus have an inclusion
$$\widetilde{\mathcal{E}(K)}^{cf} \subset 
\widetilde{\mathcal{E}(K)}$$
of strict $\s$-categories over $M \times \Omega$. 

The operad $\OO$ in $M$ gives rise to another fibered category
$$\widetilde{\OO} \longrightarrow M\times \Omega,$$
whose objects are triples $(x,T,i)$, with $x\in M$, 
$T\in \Omega$ and $i : x \rightarrow \OO(T)$ is a morphism in $M$.  \\

\textbf{Step 2:} We consider the strict $\s$-functors (i.e. simplicial 
functors) over $M\times \Omega$ 
$$\phi : \widetilde{\OO} \longrightarrow \widetilde{\mathcal{E}(K)}.$$
These functors form a category where the morphisms are the simplicial
natural transformations over $M\times \Omega$. Among these we consider the
strict $\s$-functors $\phi$ which are \emph{(strictly) cartesian in the $M$ direction}. By 
definition these are the simplicial functors $\phi$ with the following property: for 
all $(x,T,i)$ in $\widetilde{\OO}$, the canonical morphism 
$(x,T,i) \rightarrow (\OO(T),T,id_{\widetilde{\OO}(T)})$ given by $i$ is sent
to an isomorphism in $\widetilde{\mathcal{E}(K)}$. 

We denote by $\mathcal{M}$ the category of all strict $\s$-functors over $M\times \Omega$,
$\phi : \widetilde{\OO} \longrightarrow  \widetilde{\mathcal{E}(K)}$, cartesian 
in the $M$ direction. This category is easily seen to be described, up to a canonical 
equivalence, as follows. For a morphism $f : T' \rightarrow T$ in $\Omega$, and an 
object $x \in M$,  
we have the base change functor
$$f^* : \mathcal{E}(K)(x,T) \longrightarrow \mathcal{E}(K)(x,T')$$
which is left Quillen. In particular, for $x= \OO(T)$, we have a left Quillen
functor
$$f^* : \mathcal{E}(K)(\OO(T),T) \longrightarrow \mathcal{E}(K)(\OO(T),T').$$
The morphism $f$ also induces a morphism in $M$, $f^* : \OO(T) \rightarrow \OO(T')$
and thus a restriction functor on comma model categories
$$f_{!} : M/\OO(T) \longrightarrow M/\OO(T'),$$
defined by sending $x \rightarrow \OO(T)$ to the composition $x \rightarrow \OO(T')$. This
induces another left Quillen functor
$$f_{!} : \mathcal{E}(K)(\OO(T),T)\longrightarrow \mathcal{E}(K)(\OO(T'),T).$$
We note that we have an equality of left Quillen functors
$$f_{!}f^{*}=f^{*}f_{!} : \mathcal{E}(K)(\OO(T),T)\longrightarrow \mathcal{E}(K)(\OO(T'),T'),$$
and this unique functor will be simply denoted by $f^{\times}$. The association
$T \mapsto \mathcal{E}(K)(\OO(T),T)$ and $f \mapsto f^{\times}$ defines a pseudo-functor
$\mathcal{G}$ from $\Omega^{op}$ to the 2-category of model categories and left Quillen
functors. The Grothendieck's construction provides a fibered category
$$\widetilde{\mathcal{G}} \longrightarrow \Omega.$$
The category $\mathcal{M}$ can then be naturally identified with the category of 
(eventually non cartesian) sections of the projection $\widetilde{\mathcal{G}} \longrightarrow \Omega$.
To put things differently, an object of $\mathcal{M}$ consists of the following data.

\begin{enumerate}

\item For a tree $T \in\Omega$ an object $x_{T} \in \mathcal{E}(K)(\OO(T),T)$,

\item For a morphism of trees $f : T' \rightarrow T$ a morphism in $\mathcal{E}(K)(\OO(T'),T')$
$$\alpha_{f} : f^{\times}(x_{T}) \longrightarrow x_{T'}.$$

\end{enumerate}

The data $\{x_{T},\alpha_{f}\}$ is moreover required to satisfy the usual cocycle condition
with respect to composition and identity in $\Omega$. We note that by adjunction the morphism
$\alpha_{f}$ might also be described as a morphism in $\mathcal{E}(K)(\OO(T),T')$
$$f^{*}(x_{T}) \longrightarrow x_{T'} \times_{\OO(T')}\OO(T).$$
A morphism betwee two such data, $\{x_{T},\alpha_{f}\}$ and $\{x_{T},\beta_{f}\}$, is simply 
the data of a family of morphisms $x_{T} \longrightarrow y_{T}$ that commutes with the
$\alpha_{f}$'s and $\beta_{f}$'s.

It is a standard fact that such categories of sections $\mathcal{M}$ are endowed with
natural model structure for which the fibrations and the equivalences are 
defined on the underlying objects $x_{T}$. The general result is subsummed in the following
lemma. 

\begin{lem}\label{lsect}
Let $I$ be a small category and 
$$\mathcal{G} : I^{op} \longrightarrow Cat$$
a pseudo-functor. We assume that for all $i\in I$ the category 
$\mathcal{G}(i)$ is endowed with a combinatorial and simplicial model structure, and 
that for all morphism $u : i\rightarrow j$ the base change functor
$u^{*} : \mathcal{G}(j) \longrightarrow \mathcal{G}(i)$ is a simplicial left Quillen functor.
Then, there exists a model structure on the category $\mathcal{M}$ of sections
of the Grothendieck's construction $\int_{I}\mathcal{G} \longrightarrow I$, for which
the fibrations and the equivalences are defined levelwise in the model categories
$\mathcal{G}(i)$.
\end{lem}

This lemma is a slight generalisation of the existence of the projective model 
structure on categories of diagrams in combinatorial model categories, which corresponds
to the case where the pseudo-functor $\mathcal{G}$ is constant. The general case 
follows the same argument and is left to the reader (see also \cite[\S 17]{hisi}). \\

We thus have a model category $\mathcal{M}$ of strict $\s$-functors 
$\phi : \widetilde{\OO} \longrightarrow  \widetilde{\mathcal{E}(K)}$ over $M\times \Omega$
and cartesian 
in the $M$ direction. We note that for all tree $T\in \Omega$ there is a forgetful functor
$$\mathcal{M} \longrightarrow \mathcal{E}(K)(\OO(T),T)$$
which is right Quillen by definition. This functor also preserves 
cofibrations and in particular cofibrant objects, as this can be checked directly 
using the small object argument. In particular, a cofibrant and fibrant object in $\mathcal{M}$
possesses cofibrant and fibrant components in each $\mathcal{E}(K)(\OO(T),T)$. \\

\textbf{Step 3:} We now produce a specific object
$a \in \mathcal{M}$ as follows. We start by a tree of the form 
$T(n)$ (and we note as usual $\OO(n)$ for $\OO(T(n))$). 
The object $\OO(n+1)$ naturally lies over $\OO(n)$
by \emph{forgeting the last leaf}: this is the composition morphism
$$\OO(n+1) \simeq \OO(n+1) \times \OO(0) \longrightarrow \OO(n),$$
which on the level of trees correspond to gluing $T(0)$ on the last leaf of 
$T(n+1)$ and then contracting the corresponding inner edge. This defines
$\OO(n+1)$ as an object in $M/\OO(n)$. Moreover,
the object $\OO(n+1)$ also appears in a diagram in $M/\OO(n)$
$$\xymatrix{
\displaystyle{\coprod_{n}}\OO(n) \times \OO(2) \ar[r] \ar[rd] & \OO(n+1) \ar[d] & \ar[ld] \ar[l] \OO(2) \times 
\OO(n) \\
 & \OO(n), &}$$
where the vertical maps are natural projections and horizontal maps are given by 
the possible composition 
as explained right after the statement of theorem \ref{t1}. We get this way the desired
object $a_{T(n)} \in \mathcal{E}(K)(\OO(n),T(n))$.

For a general tree $T \in \Omega$, we have two product decompositions
$$\OO(T) = \prod_{v\in V(T)} \OO(T(v)) \qquad \mathcal{E}(K)(\OO(T),T)= \prod_{v \in V(T)}\mathcal{E}(K)
(\OO(T),T(v)),$$
where $v$ runs over the vertices of $T$ and $T(v)$ is the sub-tree of $T$ having $v$ as a unique
vertex and all edges adjacents to $v$ (the star in $T$ at $v$). In particular, the projections
$\OO(T) \rightarrow \OO(T(v))$ induce pull-back functors
$\mathcal{E}(K)(\OO(T(v)),T(v)) \rightarrow \mathcal{E}(K)(\OO(T),T(v))$. By considering 
the direct product of these
functors we have an external product functor
$$\boxtimes_{v} : \prod_{v\in V(T)}\mathcal{E}(K)(\OO(T(v)),T(v)) \longrightarrow \prod_{v \in V(T)}
\mathcal{E}(K)(\OO(T),T(v)) = \mathcal{E}(K)(\OO(T),T).$$
By definition, the object $a_{T}$ is defined by the formula
$$a_{T}:=\boxtimes_{v}(\{a_{T(v)}\}) \in \mathcal{E}(K)(\OO(T),T).$$

It remains to define the transition morphism $\alpha_{f}$ for $f : T' \rightarrow T$ in $\Omega$. 
For this we focus on several types of morphisms which all together generate
morphisms in $\Omega$. These different types of morphisms have already been recalled and are: 
automorphisms, Segal (or inner) faces, outer faces and units. Moreover, we have by definiton
the two product decompositions
$$\OO(T) = \prod_{v\in V(T)} \OO(T(v)) \qquad \mathcal{E}(K)(\OO(T),T)= \prod_{v \in V(T)}\mathcal{E}(K)
(\OO(T),T(v)),$$
and we will define the $\alpha_{f}$ compatibly. It will therefore
be enough to deal with the case of morphisms between very simple trees, cases that we will treat below. \\

\textbf{Outer faces:} Let $f : T' \rightarrow T$ be an outer face. We have
$\OO(T)\simeq \OO(T') \times \OO(T'')$ and a splitting $V(T)\simeq V(T') \coprod V(T'')$
(the tree $T'$ is the maximal subtree of $T$ with vertices in $V(T') \subset V(T)$, 
as well as for $T''$).
The corresponding base change functor
$$f^{*} : \mathcal{E}(K)(\OO(T),T)= \prod_{v \in V(T)}\mathcal{E}(K)
(\OO(T),T(v)) \longrightarrow \mathcal{E}(K)(\OO(T'),T')= \prod_{v \in V(T')}\mathcal{E}(K)
(\OO(T'),T'(v))$$
is obtained as the retriction along the projection $\OO(T) \rightarrow \OO(T')$ and
the projection onto the factors corresponding to the inclusion
$V(T') \subset V(T)$. The
natural morphism $\alpha_{f} : f^{*}(a_{T})\longrightarrow a_{T'}$ is then simply adjoint 
of the canonical isomorphism in $\mathcal{E}(K)(\OO(T),T)$ (which exists by definition)
$$\alpha_{f} : f^{*}(a_{T}) \simeq a_{T'} \times_{\OO(T')} \OO(T).$$

\textbf{Automorphisms:} Let $f : T=T(n) \longrightarrow T(n)$ be an automorphism of 
the tree $T(n)$ with a single vertex and $n$ leaves. 
The automorphism $f$ simply consists of a relabelling of the $n$ leaves of $T$.
The corresponding functor
$$f^{*} : \mathcal{E}(K)(\OO(n),T(n)) \longrightarrow 
\mathcal{E}(K)(\OO(n),T(n))$$
is obtained by the natural diagonal action of the symmetric group $\Sigma_{n}$ on $\OO(n)$ and on
$\mathcal{E}(K)(x,T(n))$ (for any $x$). The object $a_{T}$ is here given by the diagram
$$\xymatrix{
\displaystyle{\coprod_{n}}\OO(n) \times \OO(2) \ar[r] \ar[rd] & \OO(n+1) \ar[d] & \ar[ld] \ar[l] \OO(2) \times 
\OO(n) \\
 & \OO(n), &}$$
which is naturally $\Sigma_{n}$-equivariant, where $\Sigma_{n}$ acts simultanously on 
$\OO(n)$, $\OO(n+1)$ and on the factors of the coproduct $\coprod_{n}$. This implies that 
$a_{T}$ is naturally endowed with a structure of an $\Sigma_{n}$-object in 
$\mathcal{E}(K)(\OO(n),T(n))$, and in particular for a given automorphism $f$ we have a natural 
isomorphism
$$\alpha_{f} : f^{*}(a_{T}) \simeq a_{T}.$$

\textbf{Units:} Let $f : T(1) \longrightarrow \eta$ be the unit map in $\Omega$ (recall that 
$\eta$ is the tree with no vertices). The functor 
$$f^{*} : \mathcal{E}(K)(\OO(\eta),\eta)=* \longrightarrow \mathcal{E}(K)(\OO(1),T(1))$$
is given by an object in $\mathcal{E}(K)(\OO(1),T(1))=\mathcal{E}(K)(*,T(1))$, which is the 
following diagram of identities of objects in $M$
$$\xymatrix{
\OO(2) \ar[r] \ar[rd] & \OO(2) \ar[d] & \ar[l] \ar[ld] \OO(2) \\
 & \bullet & }$$
where these are considered as objects over $\OO(1)=*$. One the other hand, the object
$a_{1}$ in $\mathcal{E}(K)(*,T(1))$ is by definition the following diagram
of objects over $\OO(1)$. 
$$\xymatrix{
\OO(2)\simeq \OO(1) \times \OO(2) \ar[r] \ar[rd] & \OO(2) \ar[d] & \ar[ld] \ar[l] \OO(2) \times 
\OO(1)\simeq \OO(2) \\
 & \OO(1). &}$$
There exists an obvious natural isomorphism
$$\alpha_{f} : f^{*}(*) \simeq a_{1}$$
in the category $\mathcal{E}(K)(\OO(1),T(1))$. \\

\textbf{Segal faces:} Let $n$ and $m$ be two integers, and $T$ be the tree
with two vertices and $n+m$ leaves obtained by gluing the tree $T(n)$ on 
the i-th leaf of $T(m)$. By relabelling the leaves we can assume that $i=1$. 
We have a Segal face $f : T(n+m-1) \longrightarrow T$ which contracts the unique
inner edge of $T$. 

We have $\OO(T)=\OO(n)\times \OO(m)$, and as well
$$\mathcal{E}(K)(\OO(T),T)=\mathcal{E}(K)(\OO(n)\times \OO(m),T(n)) \times
\mathcal{E}(K)(\OO(n)\times \OO(m),T(m)).$$

The functor $f^{*}$ is then a functor
$$\hspace{-5mm} \mathcal{E}(K)(\OO(n)\times \OO(m),T(n)) \times
\mathcal{E}(K)(\OO(n)\times \OO(m),T(m)) \longrightarrow
\mathcal{E}(K)(\OO(n)\times \OO(m),T(n)) \times
\mathcal{E}(K)(\OO(n+m-1),T(n+m-1))$$
obtained as follows. For two objects 
$$\xymatrix{
\displaystyle{\coprod_{n}}(\OO(n)\times \OO(m) \times \OO(2)) \ar[r] \ar[rd] & F \ar[d] & \ar[ld] \ar[l] 
\OO(2) \times \OO(n)\times \OO(m) \\
 & \OO(n)\times \OO(m) , &}$$
$$\xymatrix{
\displaystyle{\coprod_{m}}(\OO(n)\times \OO(m) \times \OO(2)) \ar[r] \ar[rd] & G \ar[d] & \ar[ld] \ar[l] 
\OO(2) \times \OO(n)\times \OO(m) \\
 & \OO(n)\times \OO(m) , &}$$
we first form the composition as an object in $\mathcal{E}(K)(\OO(n)\times \OO(m),T(m+n-1))$
$$\xymatrix{
\displaystyle{\coprod_{n+m-1}}(\OO(n)\times \OO(m) \times \OO(2)) 
\ar[r] \ar[rd] & H \ar[d] & \ar[ld] \ar[l] 
\OO(2) \times \OO(n)\times \OO(m) \\
 & \OO(n)\times \OO(m) , &}$$
where $H$ is the push-out 
$$H=F\coprod_{\OO(2) \times \OO(n)\times \OO(m)}G,$$
and then we consider the cospan
$$\xymatrix{
\displaystyle{\coprod_{n+m-1}}(\OO(n)\times \OO(m) \times \OO(2)) 
\ar[r] & H  & \ar[l] 
\OO(2) \times \OO(n)\times \OO(m)}$$
as living over $\OO(n+m-1)$ via the composition morphism $\OO(n) \times \OO(m) \rightarrow \OO(n+m-1)$. 

The object $f^{*}(a_{T})$ can thus be described as the following cospan
$$\xymatrix{
\displaystyle{\coprod_{n+m-1}}(\OO(n)\times \OO(m) \times \OO(2)) \ar[r] \ar[rd] & f^{*}(a_{T}) 
\ar[d] & \ar[ld] \ar[l] 
\OO(2) \times \OO(n)\times \OO(m) \\
 & \OO(n+m-1) , &}$$
where $f^{*}(a_{T})$ is defined as a push-out
$$f^{*}(a_{T})):=
\OO(n+1)\times \OO(m)   \coprod_{\OO(n)\times \OO(m) \times \OO(2)} \OO(n)\times \OO(m+1),$$
where the two morphisms 
$$\OO(n)\times \OO(m) \times \OO(2) \rightarrow \OO(n+1)\times \OO(m) \qquad 
\OO(n)\times \OO(m) \times \OO(2) \rightarrow \OO(n)\times \OO(m+1)$$
come from the composition of $\OO(2)$ with either $\OO(n)$ or $\OO(m)$. 

On the other hand, the object $a_{T(n+m-1)}$ is by definition the cospan
$$\xymatrix{
\displaystyle{\coprod_{n+m-1}}(\OO(n+m-1) \times \OO(2)) \ar[r] \ar[rd] & \OO(n+m)
\ar[d] & \ar[ld] \ar[l] 
\OO(2) \times \OO(n+m-1) \\
 & \OO(n+m-1). &}$$
There is a canonical morphism in $\mathcal{E}(K)(\OO(n+m-1),T([n+m-1]))$
$$\alpha_{f} : f^{*}(a_{T}) \longrightarrow a_{T(n+m-1)},$$
which on the middle terms is given by 
$$\OO(n+1)\times \OO(m)  \coprod_{\OO(n)\times \OO(m) \times \OO(2)} \OO(n)\times \OO(m+1)
\longrightarrow \OO(n+m)$$
obtained by using the composition morphisms
$$\OO(n+1)\times \OO(m) \longrightarrow \OO(n+m) \qquad
\OO(n)\times \OO(m+1) \longrightarrow \OO(n+m).$$

This last case finishes the definition of the transition morphism
$\alpha_{f}$, and thus of the object $a \in \mathcal{M}$ as required. \\

\textbf{Step 4:} We have constructed an object $a \in \mathcal{M}$. We consider
a fibrant and cofibrant replacement $ad$ of $a$ in $\mathcal{M}$. This is
still a strict $\s$-functor over $M\times \Omega$ 
$$ad : \widetilde{\OO} \longrightarrow \widetilde{\mathcal{E}(K)},$$
but now factors throught the full sub-$\s$-category 
$\widetilde{\mathcal{E}(K)}^{cf}$ of cofibrant and fibrant objects. We get this way 
an $\s$-functor over $M\times \Omega$
$$ad : \widetilde{\OO} \longrightarrow \widetilde{\mathcal{E}(K)}^{cf}.$$
It is not hard to check that this morphism is fibered in the $M$ direction and thus correspond to 
a very lax  morphism of $\s$-categorical operads over $\T$
$$ad : \OO \longrightarrow \underline{End}^{\sqcup}(\OO(2)).$$
Finally, by construction this very lax morphism is also checked to be a lax
morphism in the sense of definition \ref{d2} and to satisfies the conditions of the theorem \ref{t1}. 
\hfill $\Box$ \\

\begin{df}\label{d5}
Let $\OO$ be an $\s$-operad in an $\s$-topos $\T$ with $\OO(0)\simeq \OO(1)\simeq *$. The lax morphism
of categorical $\s$-operads in $\T$
$$ad : \OO \longrightarrow \underline{End}^{\sqcup}(\OO(2))$$
of theorem \ref{t1} is called the \emph{adjoint representation of $\OO$ on $\OO(2)$}. For
an object $X \in \T$, the composition
$$ad^{X} : \xymatrix{\OO \ar[r]^-{ad} &  \underline{End}^{\sqcup}(\OO(2)) \ar[r]^-{\mathbb{D}_{X}} &
\underline{End}^{\times}(X^{\OO(2)})}$$
is called the \emph{co-adjoint representation of $\OO$ on $X^{\OO(2)}$} (where $\mathbb{D}_{X}$ is the
duality morphism of Definition \ref{d4}). 
\end{df}

The adjoint and co-ajoint representations of an operad $\OO$ are functorial with respect to 
$\OO$ in a rather obvious manner. We leave the reader to construct this functoriality (which 
follows from the functoriality of the construction of the object
$a$ during the proof of theorem \ref{t1}). \\

Let $\OO$ be an $\s$-operad in an $\s$-topos $\T$ with $\OO(0)\simeq \OO(1)\simeq *$. Theorem 
\ref{t1} shows the existence of the adjoint representation of $\OO$
$$ad : \OO \longrightarrow \underline{End}^{\sqcup}(\OO(2)),$$
as well as for any object $X \in \T$ of the co-adjoint representation
$$ad^{X} : \OO \longrightarrow \underline{End}^{\times}(X^{\OO(2)}).$$
The morphisms $ad$ and $ad^{X}$ are lax morphisms of categorical $\s$-operads
over $\T$ and we will mainly be interested in the special case for which these
morphisms are genuine (non-lax) morphisms of categorical $\s$-operads. 

\begin{df}\label{d6}
Let $\OO$ be an $\s$-operad in an $\s$-topos $\T$ such that $\OO(0)\simeq \OO(1)\simeq *$.
\begin{enumerate}

\item The $\s$-operad $\OO$ is \emph{of configuration type} if the corresponding adjoint 
representation $ad$ of $\OO$ on $\OO(2)$ is a morphism of categorical $\s$-operads. 

\item Let $\mathcal{F} \subset \T$ be a class of objects in $\T$. The $\s$-operad $\OO$ is \emph{of 
configuration type with respect to $\mathcal{F}$} if for all $X \in \mathcal{F}$ 
the corresponding co-adjoint representation $ad^{X}$ of $\OO$ on $X^{\OO(2)}$ 
is a morphism of categorical $\s$-operads. 
\end{enumerate}
\end{df}

For an object $X$, the morphism $ad^{X}$ is defined to be the composition
$$ad^{X} : \xymatrix{\OO \ar[r]^-{ad} &  \underline{End}^{\sqcup}(\OO(2)) \ar[r]^-{\mathbb{D}_{X}} &
\underline{End}^{\times}(X^{\OO(2)}),}$$
where $\mathbb{D}_{X}$ is the duality morphism of definition \ref{d4}. The morphism $\mathbb{D}_{X}$
is an actual, non-lax, morphism of categorical $\s$-operads in $\T$, and therefore
we see that if $\OO$ is of configuration type then 
it is of configuration type with respect to any classes of objects in $\T$. The converse is true by 
the Yoneda lemma: $\OO$ is of configuration type if it is of configuration type with respect
to the class of all objects in $\T$. \\

The following proposition is a consequence of the definitions and of the construction
of the morphism $ad$ of theorem \ref{t1}.

\begin{prop}\label{p1}
Let $\OO$ be an $\s$-operad in $\T$ with $\OO(0)\simeq \OO(1)\simeq *$. Then 
$\OO$ is of configuration type if and only if for all pairs of integers $n\geq 2$ and $m\geq 2$
the following square of objects in $\T$
$$\xymatrix{
\displaystyle{\OO(n+1)\times \OO(m) \coprod_{\OO(2) \times \OO(n) \times \OO(m)} \OO(n)\times \OO(m+1)}
\ar[d] \ar[r] & \OO(m+n) \ar[d] \\
\OO(n)\times \OO(m) \ar[r] & \OO(m+n-1)}$$
is cartesian in the $\s$-category $\T$. 
\end{prop}

\textit{Proof:} It is a consequence of the proof of the theorem and of the construction of
the morphism $ad$. Recall that we choose a model topos $M$ as a model for $\T$, and
well as a strict operad $\OO$ in $M$ with $\OO(0)=\OO(1)=*$. The lax morphism $ad$ is defined
first by defining an object $a$ of $\mathcal{M}$: a family of objects
$a_{T} \in \mathcal{E}(\OO(T),T)$ with compatible transitions morphisms
$\alpha_{f} : f^{*}(a_{T}) \longrightarrow a_{T'}$ for all morphism of trees
$f : T' \rightarrow T$. The morphism $ad$ is defined by chosing a fibrant and replacement
model for $a$ as an object in $\mathcal{M}$. By analysing how the 
transition morphisms $\alpha_{f}$ have been defined during the proof of the theorem \ref{t1} we
see that the $\s$-functor over $M \times\Omega$
$$ad : \widetilde{\OO} \longrightarrow \widetilde{\mathcal{E}(\OO(2))^{cf}}$$
preserves almost all cartesian morphisms in $\widetilde{\OO}$. Indeed, it does preserve
all cartesian morphisms over all morphisms in $M\times \Omega$ of the form
$(u,f)$ with $u$ arbitrary and $f$ either an automorphism, and degeneracies or an outer
face. The only remaining case to check is the case of Segal faces, which precisely 
gives the condition of the proposition.  \hfill $\Box$ \\

The dual version of the proposition is the following result. Recall that a commutative
square in $\T$
$$\xymatrix{
x \ar[r] \ar[d] & y \ar[d] \\
z \ar[r] & t}$$
is cocartesian with respect to a class of objects $\mathcal{F} \subset \T$ if for all 
$X\in \mathcal{F}$ the commutative square in $\T$
$$\xymatrix{
X^{x} & X^{y} \ar[l] \\
X^{z} \ar[u] & X^{t} \ar[l] \ar[u]}$$
is cartesian. In the same way, for $S\in \T$, a commutative square in $\T/S$
$$\xymatrix{
x \ar[r] \ar[d] & y \ar[d] \\
z \ar[r] & t}$$
is cocartesian over $S$ with respect to a class of objects $\mathcal{F} \subset \T$ if for all 
$X\in \mathcal{F}$ the commutative square in $\T/S$
$$\xymatrix{
X^{x} & X^{y} \ar[l] \\
X^{z} \ar[u] & X^{t} \ar[l] \ar[u]}$$
is cartesian in $\T/S$ (where now the exponentiation is taken relative to $S$, i.e. inside
the $\s$-category $\T/S$).

\begin{prop}\label{p2}
Let $\OO$ be an $\s$-operad in $\T$ with $\OO(0)\simeq \OO(1)\simeq *$
and $\mathcal{F}\subset \T$ a class of objects in $\T$. Then 
$\OO$ is of configuration type with respect to $\mathcal{F}$ if and only 
if for all pairs of integers $n\geq 2$ and $m\geq 2$
the following square of objects in $\T/\OO(n)\times \OO(m)$
$$\xymatrix{
\OO(2) \times \OO(n) \times \OO(m) \ar[r] \ar[d] & \OO(n+1)\times \OO(m) \ar[d] \\
\OO(n)\times \OO(m+1) \ar[r] & 
\OO(n+m) \times_{\OO(n+m-1)} (\OO(n)\times \OO(m))}
$$
is cocartesian over $\OO(n)\times \OO(m)$ with respect to $\mathcal{F}$. 
\end{prop}

\textit{Proof:} It is the same observation than for proposition \ref{p1} keeping in mind the
definition of the duality morphism $\mathbb{D}_{X}$. \hfill $\Box$ \\

\section{Derived category of the space of branes}

For the homotopy theory and the $\s$-category theory of dg-categories we refer to \cite{to1,to4,to5,ro}, from which
we will use the terminology and the notations. \\

In this section we fix $k$ a commutative ring of caracteristic zero. We will
use the notion of dg-categories and dg-functors (over $k$), as well as the notion of symmetric mono\"\i dal 
dg-categories and symmetric mono\"\i dal dg-functors (over $k$) as studied in \cite{mo}\footnote{We use 
non-closed symmetric mono\"\i dal dg-categories, whose theory is a bit simpler than 
the theory exposed in \cite{mo}.}. 
The category of symmetric mono\"\i dal dg-categories will be denoted
by $dg-Cat^{\otimes}$, and the category of dg-categories by $dg-Cat$. Both categories, 
$dg-Cat$ and $dg-Cat^{\otimes}$ have model structures (see \cite{mo}). These model structures 
are such that the forgetful functor 
$$dg-Cat^{\otimes} \longrightarrow dg-Cat$$
is a right Quillen functor.

The $\s$-categories associated to the model categories $dg-Cat$ and $dg-Cat^{\otimes}$
will be denoted by 
$$L(dg-Cat)=:\dgcat \qquad L(dg-Cat^{\otimes})=:\mdgcat.$$
From now the expressions \emph{dg-functor} and \emph{symmetric
mono\"\i dal dg-functors} will refer to morphisms in these $\s$-categories (unless
we use some specific expression such as \emph{strict dg-functor}).

We remind from \cite{to1} the notion of \emph{locally presentable}, or simply \emph{presentable},
dg-categories, which by definition are the dg-categories quasi-equivalent to left Bousfield
localisation of dg-categories of dg-modules over a small dg-category. A dg-functor 
between two such presentable dg-categories $T \rightarrow T'$ is \emph{continuous} if it 
preserves (possibly infinite) sums. Presentable dg-categories and continuous dg-functors define
a non-full sub-$\s$-category of $\dgcat$
$$\dgcat_{pr} \subset \dgcat.$$

A symmetric mono\"\i dal dg-category $A$ is \emph{presentable} if its underlying dg-category is so, and
if furthermore the dg-functor defining the mono\"\i dal structure
$$\otimes : A \otimes A \longrightarrow A$$
is continuous independently in each variable. By definition, a symmetric mono\"\i dal dg-functor between two 
presentable symmetric mono\"\i dal
dg-categories $A \rightarrow B$ is \emph{continuous} if and only if it is so 
when considered as a dg-functor between presentable dg-categories. Presentable symmetric
mono\"\i dal dg-categories and continuous symmetric mono\"\i dal dg-functors define 
a non-full sub-$\s$-category of $\mdgcat$
$$\mdgcat_{pr} \subset \mdgcat.$$ 

By definition, if $A$ is presentable symmetric mono\"\i dal dg-category, then 
the dg-functor $\otimes : A\otimes A \rightarrow A$
descends to a continuous dg-functor
$$\otimes : A\hat{\otimes} A \longrightarrow A$$
where $\hat{\otimes}$ denotes the tensor product of presentable dg-categories defined in \cite{to1}. 
Recall that for two presentable dg-categories $A$ and $B$, the dg-category
$A\hat{\otimes} B$ is presentable and comes equiped with a dg-functor
$$A\otimes B \longrightarrow A\hat{\otimes} B$$
which commutes with sums individually in each variable (we will use the
expression \emph{bicontinuous} or \emph{multicontinuous}), and which is universal, 
in the $\s$-category $\mdgcat$,
for this property. More precisely, for any dg-category $C$ which is cocomplete (i.e. is
triangulated and possesses sums), the restriction morphism
on the mapping spaces of dg-categories
$$Map(A\hat{\otimes} B,C) \longrightarrow Map(A\otimes B,C)$$
induces an equivalence on the subsimplicial set of continuous dg-functors on the
left hand side, and of bicontinuous dg-functor on the right hand side. A direct consequence of this
universal property is the following statement: for $A$ and $B$ two presentable dg-categories, 
the Yoneda embedding induces a fully faithful dg-functor
$$A\hat{\otimes} B \hookrightarrow \widehat{A\otimes B},$$
where $\widehat{A\otimes B}$ is the dg-category of cofibrant $A^{op}\otimes B^{op}$-dg-modules
(and is also equivalent to $\mathbb{R}\h(A^{op}\otimes B^{op},\hat{k})$). The essential image of 
this full embedding consists of all $A^{op}\otimes B^{op}$-dg-modules sending 
limits individually in each factor to limits (i.e. corresponding to
bicontinuous dg-functors $A\otimes B \rightarrow \hat{k}^{op}$). \\

We now let $\T$ be an $\s$-topos (t-complete as usual), and let 
$$\D : \T^{op} \longrightarrow \mdgcat_{pr}$$ 
be a stack of symmetric mono\"\i dal dg-categories.
By the adjoint functor theorem (see \cite{tove2}) we know that, after forgetting the mono\"\i dal structures, 
all the dg-functors $u^{*} : \D(x) \rightarrow \D(y)$
have right adjoints $u_{*} : \D(y) \rightarrow \D(x)$, but which might not be themselves
continuous. We have external product dg-functors
$$\D(x) \otimes \D(y) \longrightarrow \D(x\times y),$$
induced by the symmetric mono\"\i dal structures on $\D$, defined by the tensor product of the
two pull-backs along the projections
$$\D(x) \longrightarrow \D(x\times y) \qquad \D(y) \longrightarrow \D(x\times y).$$
It does factor as a dg-functor
$$\D(x) \hat{\otimes}\D(y) \longrightarrow \D(x\times y),$$
also refered to as the external product dg-functor.

We let $A$ be a presentable dg-category over $k$, which we assume to be locally cofibrant over $k$ 
(i.e. $A(a,b)$ is a cofibrant complex of $k$-modules for all pair of objects
$(a,b)$). We will define a categorical $\s$-operad $\underline{End}^{\otimes}(A)$
over $\T$, encoding families of multilinear endofunctors of $A$ parametrised by 
$\D$. 
We choose a presentation $\T\simeq LM$ for a model topos $M$. By the rectification for diagrams
we can represent, up to an equivalence, the $\s$-functor $\D : \T^{op} \longrightarrow \mdgcat$ by 
a functor of categories
$$\D : M^{op} \longrightarrow dg-Cat^{\otimes},$$
where $dg-Cat^{\otimes}$ is the model category of symmetric mono\"\i dal 
dg-categories. By composing this functor by a cofibrant replacement functor in $dg-Cat^{\otimes}$
we can even assume that for all $x\in M$ the dg-category $\D(x)$ has cofibrant hom complexes. 

For $x$ be an object of $M$ we define an Op-category
$\mathcal{E}(A)(x,-)$ as follows. For a finite set $I$ the category
$\mathcal{E}(A)(x,I)$ is the category of all 
$A^{\otimes I} \otimes A^{op} \otimes D(x)^{op}$-dg-modules $F$
$$\mathcal{E}(A)(x,I):=(A^{\otimes I}\otimes A^{op} \otimes \D(x)^{op})-Mod.$$
The identity 
object $\mathbf{1} \in \mathcal{E}(A)(x,*)$ is the diagonal bi-dg-module base changed
to $\D(x)^{op}$
$$\mathbf{1}:=A \otimes \D(x)^{op} \in (A\otimes A^{op} \otimes \D(x)^{op})-Mod,$$
where $A$ denotes the usual identity $A\otimes A^{op}$-dg-module: $(a,b) \mapsto A(b,a)$.
The symmetric group $\Sigma_{I}$ acts on $\mathcal{E}(x,I)$ by permutation of the 
factors of $A^{\otimes I}$. 

For a map of finite sets $u  : I \longrightarrow J$ the composition functor
$$\mu_{u} : 
\mathcal{E}(A)(x,J)\times \prod_{j\in J}\mathcal{E}(A)(x,u^{-1}(j)) \longrightarrow \mathcal{E}(A)(x,I)$$
is defined as follows. 
To an object $F \in \mathcal{E}(A)(x,J)$ and 
a family of objects $F_{j} \in \mathcal{E}(A)(x,u^{-1}(j))$, we form two dg-modules
$$\boxtimes_{j\in J}F_{j} \in (A^{\otimes I} \otimes (A^{op})^{\otimes J}\otimes \D(x)^{op})-Mod$$
$$F \in (A^{\otimes J}\otimes A^{op}\otimes \D(x)^{op})-Mod.$$
The first of these dg-modules is the external tensor product of the $F_{j}$'s, which
by definition is $\otimes_{j\in J}F_{j}$, considered as dg-module over the dg-category
$$\otimes_{j\in J}(A^{\otimes u^{j}}\otimes A^{op}\otimes \D(x)^{op}) \simeq 
A^{\otimes I}\otimes (A^{op})^{\otimes J} \otimes (\D(x)^{op})^{\otimes J},$$
and base changed to $A^{\otimes I}\otimes (A^{op})^{\otimes J} \otimes \D(x)^{op}$ by the
dg-functor 
$$\otimes^{J} : \D(x)^{\otimes J} \longrightarrow \D(x),$$
given by the mono\"\i dal structure on $\D(x)$. We consider
$F$ as left $A^{\otimes J}$-dg-module and a right $A\otimes \D(x)$-dg-module, and
$\boxtimes_{j\in J}F_{j}$ as a left $A^{\otimes I}\otimes \D(x)^{op}$-dg-module and
a right $A^{\otimes J}$-dg-module. 

We can compose these two dg-modules to get an new object
$$(\boxtimes_{j\in J}F_{j})\otimes_{A^{\otimes J}}F \in 
(A^{\otimes I}\otimes \D(x)^{op}\otimes A^{op} \otimes \D(x)^{op})-Mod.$$ 
Finally, we base change this by the dg-functor $\otimes : \D(x) \otimes \D(x) \longrightarrow \D(x)$
to get the desired object 
$$\mu_{u}(F,\{F_{j}\}) \in \mathcal{E}(A)(x,I).$$
This defines the functor $\mu_{u}$. 
We leave to the reader the definition of the natural isomorphisms (the $\beta$'s, $\alpha$'s and $\phi$'s) ,
given by the rather obvious choices,
making the data above into a categorical operad $\mathcal{E}(A)(x,-)$ in the sense of \S 2. \\

The categorical operad $\mathcal{E}(A)(x,-)$ is functorial in $x\in M$ as follows. For a morphism
$u : x\longrightarrow y$, we have a symmetric mono\"\i dal dg-functor 
$u^{*} : \D(y)^{op} \longrightarrow \D(x)^{op}$. 
This dg-functor induces base change functors on the level dg-modules and thus a functor
$$\mathcal{E}(A)(y,I) \longrightarrow \mathcal{E}(A)(x,I)$$
for various finite sets $I$, and promotes to an Op-functor
$\mathcal{E}(A)(y,-) \longrightarrow \mathcal{E}(A)(x,-)$. This defines 
a pseudo-functor $x \mapsto \mathcal{E}(A)(x,-)$ from $M^{op}$ to the 2-category of
categorical operads and Op-functors. 

Each category $\mathcal{E}(A)(x,I)$ comes equiped with a model structure,
for which the fibrations are epimorphisms and equivalences are quasi-isomorphisms of dg-modules.
This model category is not a simplicial model category. To overcome this difficulty
we consider $cs\mathcal{E}(A)(x,I)$ the category of co-simplicial objects in 
$\mathcal{E}(A)(x,I)$. This category comes equiped with a natural simplicial enrichement simply 
because it is the category of co-simplicial objects in a presentable category. 
We endow moreover the category $cs\mathcal{E}(A)(x,I)$ with a model structure as follows. 
We start with the projective levelwise model structure on $\mathcal{E}(A)(x,I)^{\Delta}$, 
the category of co-simplicial objects in $\mathcal{E}(A)(x,I)$, for which fibrations
and equivalences are defined levelwise. We perform a left Bousfield localization along all
the morphisms of $\Delta$ so that the fibrant objects are the objects
$E^{*}$ such that for all 
$u : [m] \rightarrow [n]$ in $\Delta$, the transition morphism $E^{m} \rightarrow E^{n}$ is an equivalence
in $\mathcal{E}(A)(x,I)$. With this model structure, there is 
a left Quillen functor
$$\begin{array}{ccc}
\mathcal{E}(A)(x,I) & \longrightarrow & cs\mathcal{E}(A)(x,I)\\
E  & \mapsto & \Delta E,
\end{array}$$
sending a dg-module $E$ to $C_{*}(\Delta^{*})\otimes E$, the co-simplicial dg-module
obtained by tensoring $E$ with the chains on the standard simplicial simplices. This
left Quillen functor is a Quillen equivalence whose right adjoint is the totalisation 
functor of co-simplicial dg-modules. In the sequel we will often use the left Quillen 
functor $E \mapsto \Delta E$ implicitely and will therefore allow ourself to 
consider objects in  $\mathcal{E}(A)(x,I)$ also as objects in $cs\mathcal{E}(A)(x,I)$. 

The categorical operad structure on $\mathcal{E}(A)(x,-)$ produces an induced categorical operad structure
on $cs\mathcal{E}(A)(x,-)$, and the pseudo-functor
$$\mathcal{E}(A) : M^{op} \times \Omega^{op} \longrightarrow \Opcat$$
extends to a pseudo-functor on co-simplicial objects
$$cs\mathcal{E}(A) : M^{op} \times \Omega^{op} \longrightarrow \Opcat.$$
This last pseudo-functor $cs\mathcal{E}(A)$ has now a natural simplicial enrichement, 
and its Grothendieck's construction provides 
a (strict) fibered $\s$-category
$$\widetilde{cs\mathcal{E}(A)} \longrightarrow M \times \Omega,$$
whose fiber at $(x,T) \in M\times \Omega$ is the strict $\s$-category
$$cs\mathcal{E}(A)(x,T):=\prod_{v\in V(T)}cs\mathcal{E}(A)(x,T(v)).$$
We denote by 
$\widetilde{cs\mathcal{E}(A)}^{cf}$ the full sub-$\s$-category of
$\widetilde{cs\mathcal{E}(A)}$ consisting of objects $(x,T,c)$, where $c$ is a cofibrant and fibrant 
object in $cs\mathcal{E}(A)(x,T)$. Using once again the sublemma \ref{sl1}
we see that $\widetilde{cs\mathcal{E}(A)}^{cf}$ remains fibered over $M\times \Omega$. 
Note that the fiber at $(x,T)$ of $\widetilde{cs\mathcal{E}(A)}^{cf}$ is an $\s$-category
naturally equivalent to the $\s$-category $L(\mathcal{E}(A)(x,T))$ of objects
in the model category $\mathcal{E}(A)(x,T)$. 

We denote by 
$\widetilde{cs\mathcal{E}(A)}_{r}^{cf}$ the full sub-$\s$-category of \emph{representable objects}
defined as follows.
Let $n$ be an integer and $(x,T(n),c)$ be an object
of $\widetilde{cs\mathcal{E}(A)}^{cf}$. By definition $c$ corresponds to a dg-module
over $A^{\otimes n}\otimes A^{op}\otimes \D(x)^{op}$, well defined
up to quasi-isomorphism. We say that $c$ is \emph{representable} if it satisfies the following two conditions.
\begin{enumerate}

\item The dg-module $c$ sends limits in each individual factor
of $A^{op}\otimes \D(x)^{op}$ to limits. In other words,  for all object 
$a \in A^{\otimes n}$, the $A^{op}\otimes \D(x)^{op}$-dg-module $c(a,-)$
belongs to the essential image of the full embedding
of dg-categories $A \hat{\otimes} \D(x) \hookrightarrow \widehat{A\otimes \D(x)}$. 

\item The corresponding dg-functor induced by the property above
$$c : A^{\otimes n} \longrightarrow A \hat{\otimes} \D(x)$$
is multicontinuous. 
\end{enumerate}

To put things differently, the representable objects are the triples $(x,T(n),c)$, where 
$c$ corresponds to the dg-module induced from a continuous dg-functor
$$A^{\hat{\otimes} n} \longrightarrow A \hat{\otimes} \D(x).$$

The $\s$-category $\widetilde{\mathcal{E}(A)}_{r}^{cf}$ remains 
a fibered $\s$-category over $M \times \Omega$.
This fibered $\s$-category provides an $\s$-functor
$$M^{op}\times \Omega^{op} \longrightarrow \scat,$$
whose value at $(x,T)$ is an $\s$-category equivalent to the simplicially enriched
category of cofibrant and representable objects in $\mathcal{E}(A)(x,T)$. This $\s$-categories is also
equivalent to the $\s$-category 
$$\prod_{v\in V(T)}\mathbb{R}\h_{c} (A^{\hat{\otimes} n_{v}},A\hat{\otimes}\D(x)),$$
of families of continuous dg-functors from the $A^{\hat{\otimes} n_{v}}$'s to 
$A\hat{\otimes}\D(x)$. For a morphism $u : x\rightarrow y$ in $M$, the base change $\s$-functor
along the morphism $(u,id) : (x,T) \rightarrow (y,T)$ is given by 
composition with $u^{*} : \D(y) \longrightarrow \D(x)$ on the right. In a similar fashion, the
base change along a morphism $(id,f) : (x,T) \rightarrow (x,T')$ is induced
from the compositions on the left. 

It is easy to see that 
the $\s$-functor $M^{op}\times \Omega^{op} \longrightarrow \scat$ obtained above
sends equivalences in $M$ to equivalences of $\s$-categories, and thus
defines an $\s$-functor
$$\underline{End}^{\otimes}(A) : \T^{op} \times \Omega \longrightarrow \scat.$$
By definition it does satisfy the Segal core condition in the second variable and thus
defines an $\s$-functor
$$\underline{End}^{\otimes}(A) : \T^{op} \longrightarrow \s-OpCat.$$

A priori this $\s$-functor is not a stack on $\T$, it does not send colimits in $\T$ 
to limits of $\s$-categories (because $A\hat{\otimes} -$ does not preserves limits
of presentable dg-categories in general). To solve this problem we introduce the following
notion. 

\begin{df}\label{d7}
Let $\D : \T^{op} \longrightarrow \dgcat_{pr}$ be a stack of presentable 
dg-categories and $A$ a presentable dg-category. The dg-category $A$ \emph{has
descent with respect to $\D$} if the $\s$-functor
$$A\hat{\otimes} \D(-) : \T^{op} \longrightarrow \dgcat_{pr}$$
is a stack. 
\end{df}

It is a formal consequence of duality that if $A$ is dualizable, as an object
in the symmetric mono\"\i dal 
$\s$-category of presentable dg-categories (see \cite{tove2}), then $A$ always has descent
with respect to any $\D$. For instance, this is the case when $A$ is compactly generated, 
as compactly generated dg-categories are dualizable objects (see \cite{to1,tove2}). \\

We assume from now that $A$ has descent with respect to $\D$ in the sense above. Then, 
the $\s$-functor we have constructed above
$$\underline{End}^{\otimes}(A) : \T^{op} \longrightarrow \s-OpCat,$$
is a stack in categorical $\s$-operads. 

\begin{df}\label{d8}
Let $A$ be a locally presentable dg-category with descent with respect to $\D$. 
The \emph{categorical $\s$-operad of endo-dg-functors
of $A$} is $\underline{End}^{\otimes}(A)$ defined above.
\end{df}

The categorical $\s$-operad $\underline{End}^{\otimes}(A)$ can be understood as an $\s$-functor
$$\T^{op} \times \Omega^{op} \longrightarrow \scat.$$
By construction, the value at $(x,T)$ of this $\s$-functor is given, up to equivalences of $\s$-categories,
by the following formula 
$$\underline{End}^{\otimes}(A)(x,T)\simeq \prod_{v\in V(T)}
\mathbb{R}\underline{Hom}_{c}(A^{\hat{\otimes} n_{v}},A\hat{\otimes}\D(x)).$$

For the next result we need to recall the notion of base change. Suppose that 
$$\xymatrix{
x \ar[r]^-{f} \ar[d]_-{p} & y \ar[d]^-{q} \\
z \ar[r]_-{g} & t}$$
is a commutative square of objects in $\T$. We say that the base change formula holds
(with respect to $\D$) if the natural transformation of dg-functor (i.e. 
the morphism in $\mathbb{R}\h_{c}(\D(y),\D(z))$)
$$h : g^{*}q_{*} \Rightarrow p_{*}f^{*}$$
is an equivalence (the natural transformation $h$ above is of course the natural one 
induced by adjunction between pull-backs and push-forward).

Using this notion, we will say that a morphism $g : z \rightarrow t$ in $\T$
\emph{satisfies the base change formula} (with respect to $\D$), if 
for all morphism $q : y \rightarrow t$ the pull-back square
$$\xymatrix{
x \ar[r]^-{f} \ar[d]_-{p} & y \ar[d]^-{q} \\
z \ar[r]_-{g} & t}$$
does so. \\

We now assume that $\OO$ is a categorical $\s$-operad over $\T$ with 
$\OO(0)=\OO(1)=*$, and $\mathcal{F}\subset \T$ is a class of objects
in $\T$. We suppose that the following conditions are satisfied.

\begin{enumerate}

\item The sub-$\s$-category $\mathcal{F}$ is a closed by finite limits and is a 
generating site for $\T$.

\item For all $X\in \mathcal{F}$, the dg-category $\D(X)$ is compactly generated, and for
all the morphisms $f : Y \rightarrow X$ in $\mathcal{F}$ the dg-functor
$f^{*} : \D(X) \longrightarrow \D(Y)$ preserves compact objects.

\item For any two objects $X$ and $Y$ in $\mathcal{F}$, the external product dg-functor
$$\D(X) \hat{\otimes} \D(Y) \longrightarrow \D(X\times Y)$$
is an equivalence of dg-categories.

\item All the morphisms between objects of $\mathcal{F}$ satisfies the base
change formula with respect to $\D$. 

\item The $\s$-operad $\OO$ is of configuration type with respect to 
$\mathcal{F}$.

\item For all $X \in \mathcal{F}$, and all $n$, the object
$X^{\OO(n)}$ belongs to $\mathcal{F}$. 

\item For an integer $n$, 
consider the natural projection $\OO(n+1) \rightarrow \OO(n)$
as well as the induced  morphism
$$Map_{/\OO(n)}(\OO(n+1),X) \longrightarrow \OO(n).$$
Then for all $Y\in \mathcal{F}$ and all morphism $Y \rightarrow \OO(n)$, 
the pull-back
$$Map_{/\OO(n)}(\OO(n+1),X) \times_{\OO(n)}Y \simeq Map_{/Y}(\OO(n+1)\times_{\OO(n)}Y,X)$$
belongs to $\mathcal{F}$.

\item For all $n$, the dg-category $\D(\OO(n))$ is compactly generated. 

\end{enumerate}

We start by some direct consequences of the assumptions above. 
We first note that  for every object $Z \in \T$, the natural morphism
$$\D(\OO(n)) \otimes \D(Z) \longrightarrow \D(\OO(n)\times Z)$$
induces an equivalence
$$\D(\OO(n)) \hat{\otimes} \D(Z) \simeq \D(\OO(n)\times Z).$$
Indeed, by condition $(1)$ we can write $Z$ has a colimit in $\T$ of object $Z_{i}$ in $\mathcal{F}$. We thus have
$$\D(Z) \simeq Lim \D(Z_{i})  \qquad \D(\OO(n)\times Z) \simeq Lim \D(\OO(n)\times Z_{i}).$$
By $(8)$ the dg-category $\D(\OO(n))$ is compactly generated and hence dualizable (see \cite{to1}), and thus
$\D(\OO(n))\hat{\otimes} -$ commutes with limits of presentable dg-categories. We have
$$Lim \D(\OO(n))\hat{\otimes}\D(Z_{i}) \simeq \D(\OO(n))\hat{\otimes}(Lim \D(Z_{i})) \simeq
\D(\OO(n))\hat{\otimes} \D(Z).$$
This reduces the problem to the case $Z \in \mathcal{F}$. We then argue similarly by writting
$\OO(n)$ as a colimit of objects in $\mathcal{F}$ and use the condition $(3)$ above to conclude that 
$$\D(\OO(n)) \hat{\otimes} \D(Z) \simeq \D(\OO(n)\times Z).$$

Another consequence is that for all tree $\T\in \Omega$ and every $X\in \mathcal{F}$, the
object $X^{\OO(T)}$ belongs to $\mathcal{F}$. Indeed, $\OO(T)$ is the product 
$\prod_{v\in V(T)}\OO(T(v))$, and we know that this is true for
$T(v)$ by condition $(6)$. The general case follows from the formula
$$X^{a\times b}\simeq (X^{a})^{b}$$
for two objects $a$ and $b$ in $\T$. 

\begin{cor}\label{cp4}
Under all the conditions $(1)-(8)$ above, there exists a morphism of categorical $\s$-operads over $\T$
$$\xymatrix{
B : \OO \ar[r] \ar[r] & \underline{End}^{\otimes}
(\D(X^{\OO(2)})).}$$
This morphism, evaluated at the tree $T(n)$, corresponds to the following continuous dg-functor
$$q_{*}p^{*} : \D(X^{\OO(2)})^{\hat{\otimes} n}\simeq \D((X^{\OO(2)})^{n}) \longrightarrow 
 \D(X^{\OO(2)} \times \OO(n)) \simeq \D(X^{\OO(2)})\hat{\otimes}  \D(\OO(n))$$
induced by the following span
$$\xymatrix{
(X^{\OO(2)})^{n} & Map_{/\OO(n)}(\OO(n+1),x) \ar[r]^-{q} \ar[l]_-{p} & X^{\OO(2)} \times \OO(n),}$$
induced from the diagram of theorem \ref{t1}, and where $Map_{/\OO(n)}$ denotes 
the internal hom object of $\T/\OO(n)$. 
\end{cor}

\textit{Proof:} 
We choose a model topos $M$ with an equivalence $\T \simeq LM$. We choose
a model $X\in M$ as a fibrant object. We let $Z:=X^{\OO(2)}$ the internal hom object
of morphisms from $\OO(2)$ to $Z$. The object $Z$ is fibrant. 
We choose model for $\D$ as a functor
$$\D : M^{op} \longrightarrow \mdgcat,$$
with cofibrant values. 

We recall (from the construction of the categorical $\s$-operad of spans in \S 3) the pseudo-functor 
$\mathcal{E}^{o}(Z) : M^{op} \longrightarrow \Opcat.$
In the same we let $A:=\D(Z)$ and consider the pseudo-functor
$\mathcal{E}(A) : M^{op} \longrightarrow \Opcat.$
We start by the construction of a natural transformation of pseudo-functors
$\D : \mathcal{E}^{o}(Z) \longrightarrow \mathcal{E}(A)$
as follows. For $x \in M$ and $I$ a finite set, the category 
$\mathcal{E}^{o}(Z)(x,I)$ is the opposite of the category of diagrams in $M$
$$\xymatrix{
\prod_{I} Z^{I} & \ar[l] F \ar[r] & Z\times x.}$$
We use the image by $\D$ and the external product dg-functors 
to get a cospan of dg-categories
$$\xymatrix{
\D(Z)^{\otimes I} \ar[r] & 
\D(\prod_{I} Z^{I})  \ar[r] & \D(F) & \ar[l] \D(Z\times x) & \ar[l] \D(Z)\otimes \D(x).}$$
From this we get two dg-functors
$$\xymatrix{ 
\D(Z)^{\otimes I} \ar[r]^-{f} & \D(F) & \D(Z) \otimes \D(x), \ar[l]_-{g}}$$
and thus a dg-module
$$\D(F)(g,f) \in \D(Z)^{\otimes I} \otimes \D(Z)^{op} \otimes \D(x)^{op},$$
sending a triplet of objects $(z,u,v)$ to the complex $\D(F)(g(u,v),f(z))$. 
This dg-module defines an object in $\mathcal{E}(A)(x,I)$, and this construction extends in an obvious
way to a lax Op-functor of Op-categories
$$\D_{x} : \mathcal{E}^o(Z)(x,-) \longrightarrow \mathcal{E}(A)(x,-).$$
Let $u : y \rightarrow x$ be a morphism in $M$. It is easy to construct a natural 
transformation of Op-functors
$$u^{*}\D_{x} \Rightarrow \D_{y} u^{*} : \mathcal{E}^o(Z)(x,-) \longrightarrow 
\mathcal{E}(A)(y,-).$$
All together this defines a morphism of fibered categories over $M \times \Omega$
$$\D : \widetilde{\mathcal{E}^o(Z)} \longrightarrow \widetilde{\mathcal{E}(A)}.$$
Passing to cosimplicial objects on both sides we get this way
a strict $\s$-functor between fibered strict $\s$-categories over $M\times \Omega$
$$\D : \widetilde{cs\mathcal{E}^o(Z)} \longrightarrow \widetilde{cs\mathcal{E}(A)}.$$
Taking the Grothendieck's construction we obtain another fibered strict $\s$-category
$$\mathcal{M} \longrightarrow M \times \Omega \times \Delta^{1}.$$
The pull-back over $M\otimes \Omega \times \{0\}$ is naturally equivalent to $\widetilde{cs\mathcal{E}(A)}$
whereas the pull-back over $M\otimes \Omega \times \{0\}$ gives back
$\widetilde{cs\mathcal{E}^o(Z)}$. We endow the categories
$cs\mathcal{E}^o(Z)(x,I)$ with simplicial model structures analogs to $cs\mathcal{E}(A)$, 
for which the fibrant objects are the homotopically constant diagrams of fibrant
objects. 

We let $\mathcal{M}^{cf}$ the full sub-$\s$-category of $\mathcal{M}$ consisting of
objects being fibrant and cofibrant in the fibers, as well as 
being representable in the sense given in this section during the definition 
of the $\s$-operad $\underline{End}^{\otimes}(A)$. This is a strict $\s$-category
over $M\times \Omega \times \Delta^{1}$, which is fibered in the $\Delta$ direction. 
By construction, the pull-back of $\mathcal{M}^{cf}$ over $M\times \Omega \times \{0\}$
is naturally equivalent to $\widetilde{\mathcal{E}(A)}^{cf}$, and the pull-back over
$M\times \Omega \times \{1\}$ is naturally equivalent to 
$\widetilde{\mathcal{E}^o(Z)}^{cf}$.
Moreover, the fibered $\s$-category $\mathcal{M}^{cf}_{r}$ produces
a lax morphism of categorical $\s$-operads over $\T$
$$\D : \underline{End}^{\times}(Z) \longrightarrow \underline{End}_{nc}^{\otimes}(A),$$
where $\underline{End}_{nc}^{\otimes}(A)$ is defined as in \ref{d8} but without the representable
condition (only cofibrant dg-modules, \emph{nc} stands for \emph{non-continuous}). 

We precompose the lax morphism $\D$ with the co-adjoint representation of \ref{d5}
$$ad^{X} : \OO \longrightarrow  \underline{End}^{\times}(Z),$$
in order to get a lax morphism of categorical $\s$-operads
$$B : \OO \longrightarrow \underline{End}_{nc}^{\otimes}(A).$$
The condition $(1)-(8)$ expressed before the statement of the corollary implies 
that this lax morphism factors 
throught $\underline{End}^{\otimes}(A)$, the sub categorical $\s$-operad of
$\underline{End}_{nc}^{\otimes}(A)$ consisting of representable dg-modules.
These conditions also imply that this lax morphism is in fact 
a morphism of categorical $\s$-operads over $\T$. By construction it satisfies the last
condition stated in the proposition. 
\hfill $\Box$ \\

\begin{df}\label{d9}
Let $X \in \T$ and $\OO$ an $\s$-operad of configuration type relative to a class 
$X\in \mathcal{F}\subset \T$ and satisfying 
the condition $(1)-(8)$ stated before corollary \ref{cp4}. The \emph{space of $\OO$-branes on $X$} is defined 
to be
the object 
$$B_{\OO}(X):=X^{\OO(2)} \in \T.$$ 
The morphism of categorical $\s$-operads in $\T$
$$B : \OO \longrightarrow \underline{End}^{\otimes}(\D(B_{\OO}(X)))$$
is called the \emph{brane operation morphism}.
\end{df}

\section{One application: higher formality}

In this last section we propose one application of the theorem \ref{t1} as well as its
consequence corollary \ref{cp4}. \\

For both of these applications we let $k$ be a base 
field of caracteristic zero. We set $\T=\dSt_{k}$ the $\s$-topos 
of derived stacks over $k$ in the sense of
\cite{seat}. A generating site for $\T$ consists of $\dAff_{k}$ the $\s$-category of affine derived
schemes endowed with the \'etale topology. The $\s$-category $\dAff_{k}$ is equivalent to the
opposite of the $\s$-category $\ncdga$, of non-positively graded commutative dg-algebras (over $k$).

The $\s$-functor
$$\D : \T^{op} \longrightarrow \mdgcat_{pr}$$
is defined as the left Kan extension of the $\s$-functor
$$\D : \dAff_{k}^{op} \longrightarrow \mdgcat_{pr}$$
sending an affine derived scheme $A$ to the symmetric mono\"\i dal
dg-category of cofibrant $N(A)$-dg-modules (where $N(A)$ is the commutative
dg-algebra obtained by normalisation from the simplicial commutative algebra $A$, see \cite{to1} for 
details).

For this paragraph we set $\mathcal{F} \subset \dSt_{k}$ to be the class of all
derived Artin stacks of the form $[Y/G]$, where $Y$ is a quasi-projective
derived scheme of finite presentation and $G$ is a smooth linear algebraic group over $k$ acting on $Y$. Any object
$X$ in $\mathcal{F}$ has a compactly generated dg-category $\D(X)$. Moreover, 
the full sub-$\s$-category $\mathcal{F}$ of $\dSt_{k}$ is stable by finite limits. \\

\textbf{Brane cohomology:} 
Let $\OO$ be an operad in $\T$ with $\OO(0)=\OO(1)=*$, and assume that all the properties $(1)-(8)$ needed for
corollary \ref{cp4} hold. We start by the observation that, for any $X \in \mathcal{F}$, 
the Brane operation morphism
$$B : \OO \longrightarrow \underline{End}^{\otimes}(\D(B_{\OO}(X)))$$
is \emph{unital}. In other words, there exists an object $\mathbf{1} \in \D(B_{\OO}(X))$,
which is a unit for all the operations
$$\D((B_{\OO}(X))^{\hat{\otimes n}} \longrightarrow \D((B_{\OO}(X))\hat{\otimes} \D(\OO(n)).$$
This follows easily from the fact that $\OO(0)=*$, the object 
$\mathbf{1} \in \D((B_{\OO}(X))$ being nothing else than the image of $k \in \D(k)$ 
by the degre $0$ operation
$$\D(k)=\D((B_{\OO}(X))^{\hat{\otimes 0}} \longrightarrow \D((B_{\OO}(X))\hat{\otimes} \D(OO(0))\simeq
\D(B_{\OO}(X)).$$
By construction of $B$, this object is easily described as
$$\mathbf{1}=j_{*}(\mathcal{O}_{X}) \in \D(B_{\OO}(X)),$$
where $j : X \longrightarrow B_{\OO}(X)$ is the \emph{constant brane} morphism. 

We leave the reader to formalise the existence of the unit $\mathbf{1}$ further
and to construct a slight modification of the categorical 
$\s$-operad $\underline{End}^{\otimes}(\D(B_{\OO}(X)))$
in the setting of dg-categories together with a fixed object. Let $\underline{End}^{\otimes}_{*}(\D(B_{\OO}(X)))$
be the this new categorical $\s$-operad over $\dSt_{k}$. Its values at $n$
and over an object $y\in \T$ is an
$\s$-category naturally equivalent to 
$$\underline{End}^{\otimes}_{*}(\D(B_{\OO}(X)))(y,n)\simeq
\mathbb{R}\h_{c}^{*}(\D(B_{\OO}(X))^{\hat{\otimes} n},\D(B_{\OO}(X))\hat{\otimes} \D(y)),$$
the $\s$-category of continuous dg-functors sending the object 
$\otimes^{n}\mathbf{1}$ to $\mathbf{1}\otimes \mathcal{O}_{y}$. There is a natural 
morphism of categorical $\s$-operads over $\dSt_{k}$
$$\underline{End}^{\otimes}_{*}(\D(B_{\OO}(X))) \longrightarrow \underline{End}^{\otimes}(\D(B_{\OO}(X)))$$
which simply forgets about the object $\mathbf{1}$. By construction the brane operation morphism $B$
factors naturally as 
$$B : \xymatrix{
\OO \ar[r]^-{B_{*}} & \underline{End}^{\otimes}_{*}(\D(B_{\OO}(X))) \ar[r] & 
\underline{End}^{\otimes}(\D(B_{\OO}(X))).}$$
Let $HH^{\OO}(X)$ be the associative dg-algebra of endomorphism of $\mathbf{1}$
$$HH^{\OO}(X):=\underline{End}_{\D(B_{\OO}(X))}(\mathbf{1}).$$
This dg-algebra is called the \emph{$\OO$-cohomology of $X$}. \\

\textbf{Operations on $\OO$-cohomology:} The corollary \ref{cp4} implies easily
that the $\s$-operad $\OO$ acts in a natural way on the dg-algebra $HH^{\OO}(X)$. This
fact can be formalized as follows. For any associative dg-algebra $B$, we define an $\s$-operad in $\dSt_{k}$
$\underline{End}^{ass}(B)$, of (associative algebra) endomorphisms of $B$. 
Its values on an integer $n$ and an affine derived scheme $y\in \dAff_{k}$ is the mapping space of 
$\OO(y)$-linear associative dg-algebras\footnote{This $\OO$ is not the same as the $\s$-operad $\OO$ but
represented functions on derived affine schemes.}
$$\underline{End}^{ass}(B)(y,n):=
Map_{\OO(y)-dg-alg}(B^{\otimes n}\otimes \OO(y),B \otimes \OO(y)).$$
The $\s$-operad $\underline{End}^{ass}(B)$ is the multiplicative analog of 
$\underline{End}(B)$, defined similarly but using mapping spaces of $\OO(y)$-dg-modules
instead of dg-algebras. There is a forgetful morphism
$$\underline{End}^{ass}(B) \longrightarrow \underline{End}(B).$$
There is an $\s$-category whose objects are complexes $E$ 
together with
a morphism of $\s$-operads 
$$\OO \longrightarrow \underline{End}(E),$$
called the $\s$-category of $\OO$-algebras. It will be denoted by 
$$\OO-\dga.$$
In the same way, there is an $\s$-category of
associative dg-algebras $B$ together with
a morphism of $\s$-operads 
$$\OO \longrightarrow \underline{End}^{ass}(B),$$
called the $\s$-category of associative algebras in $\OO$-algebras. 
It is equivalent to the $\s$-category $E_{1} \odot \OO-\dga$, 
where $\odot$ is the tensor product of $\s$-operads described in \cite{cm1,cm2} for instance. 

There is a morphism of categorical $\s$-operads in $\dSt_{k}$
$$\underline{End}^{\otimes}_{*}(\D(B_{\OO}(X))) \longrightarrow 
\underline{End}^{ass}(HH^{\OO}(X))$$
which consists of sending a dg-functor
$$\D(B_{\OO}(X))^{\hat{\otimes} n} \longrightarrow \D(B_{\OO}(X))\hat{\otimes} \D(y)$$
to the induced morphism of dg-algebras of endomorphisms of the object $\mathbf{1}^{\otimes n}$. 
The brane cohomology $HH^{\OO}(X)$ then becomes an associative dg-algebra together with 
a compatible action of $\OO$: there is an induced morphism of $\s$-operads in $\dSt_{k}$
$$\OO \longrightarrow \underline{End}^{ass}(HH^{\OO}(X)).$$
We will consider $HH^{\OO}(X)$ as an object of $E_{1}\odot \OO-\dga$.

For an equivalence $\alpha : \OO \simeq \OO'$ of $\s$-operads over $\dSt_{k}$, we have
an equivalence of $E_{1}\odot \OO$-algebras
$$HH^{\OO}(X) \simeq \alpha^{*}(HH^{\OO'}(X)),$$
where $\alpha^{*}(HH^{\OO'}(X))$ denotes $HH^{\OO'}(X)$ considered
as an $E_{1}\odot \OO$-algebra via the equivalence 
$\alpha : E_1 \odot \OO \simeq E_1 \odot \OO'$. \\

\textbf{$E_{k}$-Branes:} We now let $E_{k}$ be the topological $E_{k}$ operad, so that 
$E_{k}(n)$ is the configuration space of $n$ disjoints $k$-dimensional disks in 
the standard k-dimensional disk. We can consider $E_{k}$ as an $\s$-operad in $\dSt_{k}$, 
simply by considering spaces as constant stacks on $Spec\, k$. We observe here that 
$E_{k}$ is an $\s$-operad of configuration type. Indeed, we have to prove that 
all the diagrams
$$\xymatrix{
E_{k}(n+1) \times E_{k}(m) \displaystyle{\coprod_{E_{k}(2)\times E_{k}(n) \times E_{k}(m)}}E_{k}(n) 
\times E_{k}(m+1)
\ar[r] \ar[d]  & E_{k}(n+m) \ar[d] \\
E_{k}(n)\times E_{k}(m) \ar[r] & E_{k}(n+m-1)}$$
are homotopy cartesians (where the push-out must be understood as a homotopy push-out). For this, 
we let $(c_1,c_2)$ be a point in $E_{k}(n) \times E_{k}(m)$, with image $c \in E_{k}(n+m-1)$,
and we consider the morphism induced on the
homotopy fibers at this point. The homotopy fiber of $E_{k}(n+m) \rightarrow E_{k}(n+m-1)$
at $c$ is equivalent to the space $\mathbb{B}^{k}-\{c\}$, the complement of all the interiors of 
all the disks in the configuration $c$: it is a generalized pants, a $k$-dimensional ball 
with $n+m-1$ smaller $k$-dimensional balls removed.
On the other hand, the homotopy fiber at 
$(c_1,c_2)$ of the morphism
$$E_{k}(n+1) \times E_{k}(m) \displaystyle{\coprod_{E_{k}(2)\times E_{k}(n) \times E_{k}(m)}}E_{k}(n) 
\times E_{k}(m+1) \longrightarrow E_{k}(n)\times E_{k}(m)$$
is the homotopy push-out
$$\mathbb{B}^{k}-\{c_{1}\} \displaystyle{\coprod_{\mathbb{B}^{k}-{\{c_{0}\}}}}
\mathbb{B}^{k}-\{c_{2}\},$$
where $c_{0}$ is unique configuration of a single $k$-dimensional disk in $\mathbb{B}^{k}$.
This is  homotopy equivalent to  $\mathbb{B}^{k}-\{c\}$ by the natural map, as
$c$ is the configuration obtained by composing $c_{1}$ and $c_{2}$. 
This finishes the observation that $E_{k}$ is of configuration type.  

It is not hard to check all the properties $(1)-(8)$ of our corollary \ref{cp4}.  We thus obtain 
that $E_{k}$ acts on the dg-category $\D(B_{E_{k}}(X))$, or in other words that 
$\D(B_{E_{k}}(X))$ comes equiped with an $E_{k}$-mono\"\i dal structure. We note here
thar $E_{k}(2)$ is naturally equivalent to $S^{k-1}$, a $(k-1)$-dimensional sphere, so that 
the space of $E_{k}$-branes recovers the iterated loop spaces 
$B_{E_{k}}(X)\simeq \mathcal{L}^{k-1}X$.

On the level of $\OO$-cohomology we find that 
$$HH^{E_{k}}(X)=\underline{End}_{\D(B_{\OO}(X))}(\mathbf{1})=
\underline{End}_{\D(\mathcal{L}^{k-1}X)}(\OO_{X})$$
has natural structure of $E_{1}\odot E_{k}$-algebra. According to \cite{ha}
we have a natural equivalence of operads $E_{1}\odot E_{k} \simeq E_{k+1}$, which 
promotes $HH^{E_{k}}(X)$ into an $E_{k+1}$-algebra.  This proves the following
corollary, recovering the results of \cite{bfn} by different methods.

\begin{cor}\label{c2}
Let $k\geq 0$ be an integer, $X$ an object in $\mathcal{F}$
and $\mathcal{L}^{k-1}X=B_{E_{k}}(X)$ be the space of $E_{k}$-branes in $X$.
\begin{enumerate}
\item The dg-category $\D(\mathcal{L}^{k-1}X)$ is naturally endowed with
an $E_{k}$-mono\"\i dal structure.
\item The $E_{k}$-cohomology of $X$ 
$$HH^{E_{k}}(X)=\underline{End}_{\D(\mathcal{L}^{k-1}X)}(\OO_{X})$$
is naturally endowed with an $E_{k+1}$-algebra structure.
\end{enumerate}
\end{cor}

\textbf{$E_{k}^{uni}$-Branes:} The previous notion of $E_{k}$-branes possesses a slight modification, 
consisting of replacing spheres $S^{k-1}$ by unipotent, or \emph{affine}, spheres
$(S^{k-1}\otimes k)^{uni}$ in the sense of \cite{chaff}. For this, we remind the notion
of \emph{affine stacks}. 

Let $A$ be a non-negatively graded commutative dg-algebra. We associate to it 
a derived stack $Spec\, A$ defined as the $\s$-functor
$$Spec\, A : \xymatrix{\dAff^{op}=\ncdga \ar[r] & \Top \\
 B \ar[r] & Map_{\cdga}(A,B).}
$$
The construction $Spec\, $ defines an $\s$-functor from 
$\cdga^{\geq 0}$, of non-negatively graded commutative dg-algebras, to derived stacks. 
It can be shown that this $\s$-functor is fully faithful when restricted to 
the objects $A \in \cdga^{\geq 0}$ for which there is an augmentation 
$A \rightarrow k$ inducing an isomorphism $H^{0}(A)\simeq k$ (such dg-algebras
will be called \emph{reduced}). The essential image of this full embedding
will be called the $\s$-\emph{category of (reduced) affine stacks}.

There exists an $\s$-functor (denoted by $X \mapsto (X\otimes k)^{uni}$ in \cite{chaff})
$$\begin{array}{ccc}
\Top & \longrightarrow & \dSt_{k} \\
X & \mapsto & Spec\, (C^*(X,k))=:X^{u},
\end{array}$$
sending a space $X$ to the spectrum of its commutative dg-algebra of cochains
$C^*(X,k)$ (here we use a model for $C^{*}(X,k)$ which is a commutative
dg-algebra). 
When restricted to connected spaces, this $\s$-functor factors throught the
full sub-$\s$-category of affine stacks and is called the \emph{affinization functor}. 
The Kunneth formula implies that it preserves finite products when 
restricted to the $\s$-category $\Top_{0}^{f}$ consisting of connected spaces being furthermore 
retracts of finite CW complexes (i.e. finitely presented objects in 
the $\s$-category $\Top$). 

Let $\OO$ be a topological operad such that each space $\OO(n)$ 
is connected and weakly equivalent to a retract of a finite CW complex. We present $\OO$ has an $\s$-functor
$$\Omega^{op} \longrightarrow \Top_{0}^{f}$$
satisfying the Segal core conditions. Composing with the $\s$-functor $(-)^{u}$ we obtain another
$\s$-functor
$$\OO^{u} : \Omega^{op} \longrightarrow \dSt_{k},$$
which is now an $\s$-operad in the $\s$-topos of derived stacks. By construction we have a universal
morphism of $\s$-operads in $\dSt_{k}$ 
$$\OO \longrightarrow \OO^{u},$$
where $\OO$ is considered asa  constant $\s$-operad in $\dSt_{k}$. The following lemma follows
essentially by the universal property of the construction $X \mapsto X^{u}$ (see \cite{chaff}).

\begin{lem}\label{l3}
The morphism $\OO \longrightarrow \OO^{u}$ induces an equivalence of $\s$-categories 
of dg-algebras
$$\OO^{u}-\dga \simeq \OO-\dga.$$
\end{lem}

We use this construction to consider $E_{k}^{u}$ for $k>1$, which by definition is the 
\emph{unipotent little $k$-disks $\s$-operads}. It is an $\s$-operad of configuration
type with respect to the class $\mathcal{F}$, as this can be checked using that 
$E_{k}$ is of configuration type. Moreover, it does satisfy the conditions $(1)-(8)$ of
our corollary \ref{cp4}. For $X \in \mathcal{F}$, the space of branes with respect to $E_{k}^{u}$
is here 
$$B_{E_{k}^{u}}(X)\simeq Map(S^{k-1}_{u},X)=\mathcal{L}_{u}^{(k-1)}(X) \subset \mathcal{L}^{(k-1)}(X),$$
where $S_{u}^{k-1}=(S^{k-1})^{u}$ denotes the unipotent $(k-1)$-sphere, and
$\mathcal{L}_{u}^{(k-1)}(X)$ is the \emph{unipotent part} of the iterated 
loop space $\mathcal{L}^{(k-1)}(X)$. We thus get an $E_{k}^{u}$-action
on the dg-category $\D(\mathcal{L}_{u}^{(k-1)}(X))$. 

The $E_{k}^{u}$-cohomology is given by 
$$HH^{E_{k}^{u}}(X)=\underline{End}_{\D(\mathcal{L}_{u}^{k-1}X)}(\OO_{X})$$
and is an $E_{1}\odot E_{k}^{u}$-algebra. By construction there is an
inclusion morphism
$$HH^{E_{k}^{u}}(X) \longrightarrow HH^{E_{k}}(X)$$
corresponding to the inclusion $\mathcal{L}^{(k-1)}_{u}(X) \hookrightarrow \mathcal{L}^{(k-1)}(X)$. 
Because our assumptions insure that $X$ is a derived Artin $1$-stack
we can check that this inclusion morphism is a quasi-isomorphism because $k>1$. 
We can moreover see that this quasi-isomorphism is compatible with the
operations of $E_{k}$ and $E_{k}^{u}$ via the natural morphism $E_{k} \rightarrow E_{k}^{u}$. 
In fact more is true, this morphism of $\s$-operads induces an 
equivalence between $\s$-categories of algebras
$$E_{1}\odot E_{k}^{u}-\dga \simeq E_{1}\odot E_{k}-\dga.$$
Throught this equivalence of $\s$-categories the object 
$HH^{E_{k}^{u}}(X)$ is identified with $HH^{E_{k}}(X)$. \\

\textbf{$\mathbb{P}_{k}$-Branes:} The theory of $\mathbb{P}_{k}$-branes is obtained
by considering a formal analog of the theorie of $E_{k}^{u}$-branes. We consider the $\s$-functor
$$\begin{array}{ccc}
\Top & \longrightarrow & \dSt_{k} \\
X & \mapsto & Spec\, (H^*(X,k))=:X^{f},
\end{array}$$
which provides an $\s$-functor
$$\Top^{f}_{0} \longrightarrow \dSt_{k},$$
on the $\s$-category $\Top^{f}_{0}$ of connected and homotopically finite spaces. 
This $\s$-functor commutes with finite products and thus induces an $\s$-functor
on $\s$-operads
$$\s-Op(\Top^{f}_{0}) \longrightarrow \s-Op(\dSt_{k}),$$
denotes by $\OO \mapsto \OO^{f}$. Applied to the operad $E_{k}$ for $k>1$ we find
by definition
$$\mathbb{P}_{k}:=E_{k}^{f}.$$
The $\s$-operad is a realisation in $\dSt_{k}$ of the standard 
Poisson operads classfying Poisson dg-algebras whose brackets is of
cohomological degree $1-k$. The $\s$-category $\mathbb{P}_{k}-\dga$ can be identified
with the $\s$-category of commutative dg-algebras endowed with a Poisson 
bracket of degree $1-k$ (these are also called \emph{$k$-algebras} in \cite{gj}).

It can be checked that $\mathbb{P}_{k}$ is of configuration type with respect to 
$\mathcal{F}$ and satisfies the conditions $(1)-(8)$ of corollary \ref{cp4}. For 
$X \in \mathcal{F}$, the space of $\mathbb{P}_{k}$-branes is given by the \emph{iterated
formal loop stack}
$$B_{\mathbb{P}_{k}}(X)\simeq Map(S^{k-1}_{f},X)=\mathcal{L}^{(k-1)}_{f}X,$$
where $S^{k-1}_{f}=Spec\, H^{*}(S^{k-1})$ is the formal $(k-1)$-dimensional sphere. 
Note that $S^{k-1}_{f}$ can also be considered as 
$Spec\, k[\epsilon_{k-1}]$, where $\epsilon_{k-1}$ is a square zero function
of cohomological degree $k-1$. 

We get this way the notion of $\mathbb{P}_{k}$-cohomology of $X$
$$HH^{\mathbb{P}_{k}}(X)=\underline{End}_{\D(\mathcal{L}_{f}^{k-1}X)}(\OO_{X}),$$
which is an $E_{1}\odot \mathbb{P}_{k}$-algebra. The HKR theorem implies that 
as a complexe we have
$$HH^{\mathbb{P}_{k}}(X)\simeq \mathbb{R}\Gamma(X,Sym_{\OO_{X}}(\mathbb{T}_{X}[-k])),$$
where $\mathbb{T}_{X}$ is the tangent complex of $X$
and the $E_{1}\odot \mathbb{P}_{k}$-structure reflects the
operations encoded by brackets of vector fields on $X$. The object
$HH^{\mathbb{P}_{k}}(X)$ is a shifted and derived analog of 
polyvector fields on smooth manifold considered in \cite{ko}. \\

\textbf{The formality theorem:} We will now prove a version of the
higher formality conjecture mentioned in the introduction of \cite{ptvv}. A more
detailed argument, as well as investigations of how to weaken our hypothesis on $X$, should appear 
in a forthcioming work.

We let $X$ be a derived algebraic stack in $\mathcal{F}$. Remind that this
implies that $X$ is a quotient of a quasi-projective derived scheme
by a reductive algebraic group. Let $k>1$ and consider the two 
$\s$-operads in $\dSt_{k}$ already mentioned:
$$E_{k}^{u} \qquad \mathbb{P}_{k}=\mathbb{E}_{k}^{f}.$$
By the formality of the $E_{k}$-operad (see \cite{lv}) we can choose an equivalence of
$\s$-operads in $\dSt_{k}$
$$\alpha : \mathbb{P}_{k} \simeq E_{k}^{u}.$$
This equivalence provides an equivalence between the $\s$-category
of $E_{1}\odot E_{k}^{u}$-algebras and the 
$\s$-category of $E_{1}\odot \mathbb{P}_{k}$-algebras
$$\alpha^{*} : E_{1}\odot E_{k}^{u}-\dga \simeq E_{1}\odot \mathbb{P}_{k}-\dga.$$
As we have seen, there is a natural equivalence of $E_{1}\odot \mathbb{P}_{k}$-algebras
$$\alpha^{*}(HH^{E^{u}_{k}}(X)) \simeq HH^{\mathbb{P}_{k}}(X).$$

We have forgetful $\s$-functors to the $\s$-category of dg-Lie algebras
$$\mathbb{P}_{k}-\dga \longrightarrow \dgl \qquad E^{u}_{k}-\dga \longrightarrow \dgl$$
which are induced by the standard embedding of operads in chain complexes
$$Lie_{k} \hookrightarrow H_{*}(E_{k}) \qquad Lie_{k} \hookrightarrow C_{*}(E_{k}) .$$
On the level of underlying complexes, these forgetful $\s$-functors
send $B$ to $B[k-1]$.
The equivalence $\alpha$ automatically preserves these inclusions, and thus
the equivalence of $\s$-categories $\alpha^{*}$ fits into a commutative diagram 
of $\s$-categories
$$\xymatrix{
E_{1}\odot E_{k}^{u}-\dga \ar[rr]^-{\alpha^{*}} \ar[rd] & & E_{1}\odot \mathbb{P}_{k}-\dga \ar[dl] \\
 & E_{1}(\dgl), &
}$$
where $E_{1}(\dgl)$ is the $\s$-category of mono\"\i d objects
inside the $\s$-category of dg-Lie algebras. Here we consider
$\dgl$ as a symmetric mono\"\i dal $\s$-category with the direct product structure
as a mono\"\i dal structure: $E_{1}$-mono\"\i ds can simply be represented
as simplicial objects in $\dgl$ satisfying a Segal type condition.
We find this way an equivalence in $E_{1}(\dgl)$
$$HH^{\mathbb{P}_{k}}(X)[k-1] \simeq HH^{E_{k}}(X)[k-1].$$

\begin{lem}\label{l2}
The classifying space $\s$-functor
$$B : E_{1}(\dgl) \longrightarrow \dgl$$
is an equivalence of $\s$-categories.
\end{lem}

\textbf{Proof of the lemma:} For any pointed and presentable $\s$-category $\mathcal{C}$, 
the $\s$-category $E_{1}(\mathcal{C})$ is the full sub-$\s$-category of
$s\mathcal{C}:=Fun^{\infty}(\Delta^{op},\mathcal{C})$ consisting of simplicial
objects $M_{*}$ such that for all $n$ the Segal maps induce an equivalence 
$$M_{n} \simeq M_{1}^{n}.$$
We have an adjunction
$$B : E_{1}(\mathcal{C}) \longrightarrow \mathcal{C} : \Omega,$$
where $B$ is the colimit over $\Delta$ and 
$\Omega$ sends an object $X$ to its loop space $\Omega(X)=*\times_{X}*$. When 
$\mathcal{C}$ is stable this adjunction is automatically an equivalence of $\s$-categories. This
implies that the same is true for $\mathcal{C}=\dgl$, as this can be seen using the
forgetful $\s$-functor from dg-Lie algebras to complexes, which is conservative and
commutes with the construction $B$ and $\Omega$ (because $B$ is a sifted colimit $\s$-functor).
\hfill $\Box$ \\

The lemma implies that the both complexes
$$HH^{\mathbb{P}_{k}}(X)[k] \qquad  HH^{E_{k}}(X)[k]$$
are endowed with natural dg-Lie algebra structures, and the equivalence $\alpha$
induces an equivalence of dg-Lie algebras
$$HH^{\mathbb{P}_{k}}(X)[k] \simeq HH^{E_{k}}(X)[k].$$

To finish, by construction $HH^{E_{k}}(X)[k]$
is the shifted iterated Hochschild cohomology of $X$ endowed with 
its dg-Lie alegbra structure induced by the $E_{k+1}$-algebra structure
on $HH^{E_{k}}(X)$. On the other hand $HH^{\mathbb{P}_{k}}(X)[k]$
is the complex of shifted polyvectors
$$HH^{\mathbb{P}_{k}}(X)[k] \simeq \mathbb{R}\Gamma(X,Sym_{\OO_{X}}(\mathbb{T}_{X}[-k]))[k],$$
endowed with the dg-Lie algebra structure induced by a version of the Schouten bracket
on polyvector fields. 

\begin{cor}\label{cform}
Let $X \in \mathcal{F}$ be a derived algebraic stack and $k>1$. A choice of an 
equivalence of $\s$-operads
$$\alpha : \mathbb{P}_{k} \simeq E_{k}^{u}$$
induces an canonical equivalence of dg-Lie algebras
$$\mathbb{R}\Gamma(X,Sym_{\OO_{X}}(\mathbb{T}_{X}[-k]))[k] \simeq HH^{E_{k}}(X)[k].$$
\end{cor}

The case $k=2$ and $X$ a smooth variety, and the known relations between 
$HH^{E_{2}}$ and the deformation theory of mono\"\i dal categories (see \cite{fr}), 
answers positively the conjecture of Kapustin \cite[p. 14]{kap}.

\end{document}